# Natasha: Faster Non-Convex Stochastic Optimization Via Strongly Non-Convex Parameter

(version 5)

Zeyuan Allen-Zhu
zeyuan@csail.mit.edu
Princeton University / Institute for Advanced Study

February 2, 2017[*]

**Abstract**

Given a nonconvex function that is an average of $n$ smooth functions, we design stochastic first-order methods to find its approximate stationary points. The convergence of our new methods depends on the smallest (negative) eigenvalue $-\sigma$ of the Hessian, a parameter that describes how nonconvex the function is.

Our methods outperform known results for a range of parameter $\sigma$, and can be used to find approximate local minima. Our result implies an interesting dichotomy: there exists a threshold $\sigma_0$ so that the currently fastest methods for $\sigma > \sigma_0$ and for $\sigma < \sigma_0$ have different behaviors: the former scales with $n^{2/3}$ and the latter scales with $n^{3/4}$.

## 1 Introduction

We study the problem of composite *nonconvex* minimization:

$$\min_{x \in \mathbb{R}^d} \left\{ F(x) \stackrel{\text{def}}{=} \psi(x) + f(x) \stackrel{\text{def}}{=} \psi(x) + \frac{1}{n} \sum_{i=1}^n f_i(x) \right\} \quad (1.1)$$

where each $f_i(x)$ is *nonconvex but smooth*, and $\psi(\cdot)$ is proper convex, possibly nonsmooth. We are interested in finding a point $x$ that is an *approximate local minimum* of $F(x)$.

- The finite-sum structure $f(x) = \frac{1}{n} \sum_{i=1}^n f_i(x)$ arises prominently in large-scale machine learning tasks. In particular, when minimizing loss over a training set, each example $i$ corresponds to one loss function $f_i(\cdot)$ in the summation. This finite-sum structure allows one to perform stochastic gradient descent with respect to a random $\nabla f_i(x)$.

- The so-called *proximal* term $\psi(x)$ adds more generality to the model. For instance, if $\psi(x)$ is the indicator function of a convex set, then problem (1.1) becomes constraint minimization; if $\psi(x) = \|x\|_1$, then we can allow problem (1.1) to perform feature selection. In general, $\psi(x)$ has to be a simple function where the projection operation $\arg\min_x\{\psi(x) + \frac{1}{2\eta}\|x - x_0\|^2\}$ is efficiently computable. At a first reading of this paper, one can assume $\psi(x) \equiv 0$ for simplicity.

---

[*]V1 appeared on arXiv on this date. V2, V3, V4, V5 corrected several typos, minor mistakes, polished writing, and added citations. No new result is added since V1. After V2, we wrote a sister paper `Natasha2` [3]; it addresses online methods and is built on top of this work. For this reason, our method in this paper is renamed from `Natasha` to `Natasha1` since V3.

Many nonconvex machine learning problems fall into problem (1.1). Most notably, training *deep neural networks* and classifications with *sigmoid loss* correspond to (1.1) where neither $f_i(x)$ or $f(x)$ is convex. However, our understanding to this challenging nonconvex problem is very limited.

## 1.1 Optimization with Bounded Nonconvexity

Let $L$ be the smoothness parameter for each $f_i(x)$, meaning all the eigenvalues of $\nabla^2 f_i(x)$ lie in $[-L, L]$.[1] We denote by $\sigma \in [0, L]$ the *bounded nonconvexity* parameter of $f(x) = \frac{1}{n} \sum_{i=1}^n f_i(x)$, meaning that[2]

$$\text{all the eigenvalues of } \nabla^2 f(x) \text{ lie in } [-\sigma, L].$$

We say $f(x)$ is of $\sigma$-bounded nonconvexity (or just $\sigma$-nonconvex for short). This parameter $\sigma$ should be reminiscent of the strong-convexity parameter $\mu$ for convex optimization, where all the eigenvalues of $\nabla^2 f(x)$ lie in $[\mu, L]$ for some $\mu > 0$.

We wish to find an $\varepsilon$-approximate stationary point (a.k.a. critical point) of $F(x)$, that is

$$\text{a point } x \text{ satisfying } \|\mathcal{G}(x)\| \leq \varepsilon$$

where $\mathcal{G}(x)$ is the so-called gradient mapping of $F(x)$ (see Section 2 for a formal definition). In the special case of $\psi(\cdot) \equiv 0$, one has $\mathcal{G}(x) = \nabla f(x)$.

Since $f(\cdot)$ is of $\sigma$-bounded nonconvexity, at least when $\psi(\cdot) \equiv 0$, any $\varepsilon$-approximate stationary point is automatically also an $(\varepsilon, \sigma)$-approximate local minimum— meaning that the Hessian of the output point $\nabla^2 f(x) \succeq -\sigma \mathbf{I}$ is approximately positive semidefinite (PSD).

## 1.2 Motivations and Remarks

- We focus on optimization with bounded nonconvexity because introducing this parameter $\sigma$ allows us to perform a *more refined study* of non-convex optimization. If $\sigma$ equals $L$ then optimization with $L$-bounded nonconvexity is equivalent to the general non-convex (smooth) optimization. We hope that this encourages a new way to compare nonconvex algorithms.

- We focus only on finding stationary points as opposed to local minima, because in recent studies [1, 3, 5, 13] —see Appendix B— it is shown that finding $(\varepsilon, \delta)$-approximate local minima *reduces* to finding $\varepsilon$-approximate stationary points in functions of $O(\delta)$-bounded nonconvexity.

- Parameter $\sigma$ is often not constant and can be much smaller than $L$. For instance, second-order methods often find $(\varepsilon, \sqrt{\varepsilon})$-approximate local minima [26] and this corresponds to $\sigma = \sqrt{\varepsilon}$.

## 1.3 Known Results

Despite the widespread use of nonconvex models in machine learning and related fields, our understanding to non-convex optimization is still very limited. Until recently, nearly all research papers have been mostly focusing on either $\sigma = 0$ or $\sigma = L$:

---

[1]This definition also applies to functions $f(x)$ that are not twice differentiable, see Section 2 for details.

[2]Previous authors also refer to this notion as "approximate convex", "almost convex", "hypo-convex", "semi-convex", or "weakly-convex." We call it $\sigma$-nonconvex to stress the point that $\sigma$ can be as large as $L$ (recall any $L$-smooth function is automatically $L$-nonconvex). In our earlier versions of this paper, we have called $\sigma$ the "strong nonconvexity" parameter, but were told by some readers that it is a bad notion. We have renamed it since then, but kept the paper title unchanged.

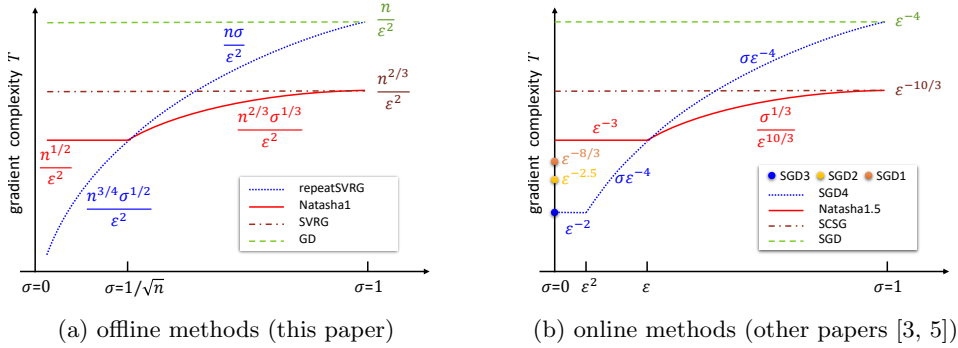

(a) offline methods (this paper)  (b) online methods (other papers [3, 5])

Figure 1: Nonconvex optimization for finding $\|\nabla f(x)\| \leq \varepsilon$ on functions with $\sigma$-bounded nonconvexity.

- If $\sigma = 0$, the accelerated SVRG method [12, 30] finds $x$ satisfying $F(x) - F(x^*) \leq \varepsilon$, in gradient complexity $\widetilde{O}\big(n + n^{3/4}\sqrt{L/\varepsilon}\big)$.[3] This result studies convex $f(x)$ and is irrelevant to this paper.
- If $\sigma = L$, the SVRG method [6] finds an $\varepsilon$-approximate stationary point of $F(x)$ in gradient complexity $O(n + n^{2/3}L/\varepsilon^2)$.
- If $\sigma = L$, full gradient descent (GD) finds an $\varepsilon$-approximate stationary point of $F(x)$ in gradient complexity $O(nL/\varepsilon^2)$.
- If $\sigma = L$, stochastic gradient descent (SGD) finds an $\varepsilon$-approximate stationary point of $F(x)$ in gradient complexity $O(L/\varepsilon^2 + L\mathcal{V}/\varepsilon^4)$ where $\mathcal{V}$ is the variance of the stochastic gradient.[4]

Throughout this paper, we refer to *gradient complexity* as the total number of stochastic gradient computations $\nabla f_i(x)$ and proximal computations $y \leftarrow \mathsf{Prox}_{\psi,\eta}(x) \stackrel{\mathrm{def}}{=} \arg\min_y \{\psi(y) + \frac{1}{2\eta}\|y - x\|^2\}$.[5]

To the best of our knowledge, even if $0 < \sigma \ll L$, it is not clear whether SGD, GD, or SVRG can take advantage of $\sigma$.[6] Very recently, it was observed by two independent groups [1, 13] —although implicitly, see Section 2.1— that for minimizing functions of $\sigma$-bounded nonconvexity, one can repeatedly regularize $F(x)$ to make it $\sigma$-strongly convex, and then apply the accelerated SVRG method to minimize this regularized function. Under mild assumption $\sigma \geq \varepsilon^2$, this approach

- finds an $\varepsilon$-approximate stationary point in gradient complexity $\widetilde{O}\big(\frac{n\sigma + n^{3/4}\sqrt{L\sigma}}{\varepsilon^2}\big)$.

We call this method repeatSVRG in this paper. Unfortunately, repeatSVRG is even slower than the vanilla SVRG for $\sigma = L$ by a factor $n^{1/3}$, see Figure 1a.

### 1.4 Our New Results

In this paper, we focus on *offline methods* which are algorithms that run in gradient complexity polynomial in $n$, but at most quadratically in $\varepsilon^{-1}$. For instance, SGD is *not* offline.

We identify an interesting dichotomy with respect to the spectrum of the nonconvexity parameter $\sigma \in [0, L]$. In particular, we showed that if $\sigma \geq L/\sqrt{n}$, then our new method `Natasha1` finds

---

[3] We use $\widetilde{O}$ to hide poly-logarithmic factors in $n, L, 1/\varepsilon$.

[4] The non-convex convergence rates of GD/SGD are not hard to prove. The rate for GD was recorded in Nesterov [25], and was perhaps first established by Polayk in 1960s. The rate for SGD first dates back to Ghadimi and Lan [19].

[5] Some authors also refer to them as incremental first-order oracle (IFO) and proximal oracle (PO) calls. In most machine learning applications, each IFO and PO call can be implemented to run in time $O(d)$ where $d$ is the dimension of the model, or even in time $O(s)$ if $s$ is the average sparsity of the data vectors.

[6] Even when $\sigma = 0$, the task of finding a point with $\|\nabla f(x)\| \leq \varepsilon$ is a non-trivial task, see [27].

an $\varepsilon$-approximate stationary point of $F(x)$ in gradient complexity

$$O\Big(n\log\frac{1}{\varepsilon} + \frac{n^{2/3}(L^2\sigma)^{1/3}}{\varepsilon^2}\Big) \ . \tag{Natasha1}$$

In other words, together with repeatSVRG, we have improved the (offline) gradient complexity for nonconvex optimization of $\sigma$-bounded nonconvexity to[7]

$$\widetilde{O}\Big(\min\Big\{\frac{n^{3/4}\sqrt{L\sigma}}{\varepsilon^2}, \frac{n^{2/3}(L^2\sigma)^{1/3}}{\varepsilon^2}\Big\}\Big) \tag{repeatSVRG}$$

and the first term in the min is smaller if $\sigma < L/\sqrt{n}$ and the second term is smaller if $\sigma > L/\sqrt{n}$. We illustrate our performance improvement in Figure 1a. Our result matches that of SVRG for $\sigma = L$, and has a simpler analysis.

## 1.5 Our Extensions

One can take a step further and ask what if each function $f_i(x)$ is $(\ell_1, \ell_2)$-smooth for parameters $\ell_1, \ell_2 \geq \sigma$, meaning that all the eigenvalues of $\nabla^2 f_i(x)$ lie in $[-\ell_2, \ell_1]$. We show that a variant of our method Natasha1$^{\text{full}}$ solves this more refined problem of (1.1) with total gradient complexity

$$O\Big(n\log\frac{1}{\varepsilon} + \frac{n^{2/3}(\ell_1\ell_2\sigma)^{1/3}}{\varepsilon^2}\Big) \tag{Natasha1$^{\text{full}}$}$$

as long as $\frac{\ell_1\ell_2}{\sigma^2} \leq n^2$. In contrast, repeatSVRG achieves (for $\sigma \geq \varepsilon^2$)

$$\widetilde{O}\Big(\frac{n\sigma + n^{1/2}((\ell_1+\ell_2)\sigma)^{1/2} + n^{3/4}(\ell_1\ell_2\sigma^2)^{1/4}}{\varepsilon^2}\Big) \tag{repeatSVRG}$$

and is worse than Natasha1$^{\text{full}}$ if $n \geq \frac{\ell_1\ell_2}{\sigma^2}$ or $n \leq \frac{\max\{\ell_1,\ell_2\}\sigma}{\min\{\ell_1^2,\ell_2^2\}}$, but better than Natasha1$^{\text{full}}$ otherwise. (We would like to point out that, in our V1-V4 of this paper, we forgot to add the term $n^{1/2}((\ell_1+\ell_2)\sigma)^{1/2}$ in the numerate for this refined statement of repeatSVRG. We have fixed it in V5.)

*Remark* 1.1. In applications, $\ell_1$ and $\ell_2$ can be of very different magnitudes. The most influential example is finding the leading eigenvector of a symmetric matrix. Using the so-called shift-and-invert reduction [17], computing the leading eigenvector reduces to the convex version of problem (1.1), where each $f_i(x)$ is $(\lambda, 1)$-smooth for $\lambda \ll 1$. Other examples include all the applications that are built on shift-and-invert, including high rank SVD/PCA [7], canonical component analysis [8], online matrix learning [9], and approximate local minima algorithms [1, 13].

**Mini-Batch Setting.** Our result generalizes trivially to the mini-batch stochastic setting, where in each iteration one computes $\nabla f_i(x)$ for $b$ random choices of index $i \in [n]$ and average them. The stated gradient complexities of Natasha1 and Natasha1$^{\text{full}}$ can be adjusted so that the factor $n^{2/3}$ is replaced with $n^{2/3}b^{1/3}$.

**Online Stochastic Setting.** In super large-scale settings, it can be desirable to design iterative algorithms with gradient complexities *independent* of $n$. Such methods are *online* methods. For instance, stochastic gradient descent (SGD) is an online method and has a convergence rate of $T = O(1/\varepsilon^4)$.

Since the original appearance of this paper, several online complexities were discovered for non-convex optimization. Lei et al. [23] proposed a variant of SVRG (and called it SCSG) with gradient complexity $O\big(\frac{1}{\varepsilon^{10/3}}\big)$. Allen-Zhu [3] generalized Natasha1 to the online setting (and called it Natasha1.5) with gradient complexity $O\big(\frac{1}{\varepsilon^3} + \frac{\sigma^{1/3}}{\varepsilon^{10/3}}\big)$. Allen-Zhu [5] proposed a variant of SGD (and called it SGD5) with gradient complexity $O\big(\frac{1}{\varepsilon^2} + \frac{\sigma}{\varepsilon^4}\big)$. We compare them in Figure 1b. These algorithms can also be upgraded for finding *approximate local minima*; for more details, see [3, 5].

---

[7]We remark here that this is under mild assumptions for $\varepsilon$ being sufficiently small. For instance, the result of [1, 13] requires $\varepsilon^2 \leq \sigma$. In our result, the term $n\log\frac{1}{\varepsilon}$ disappears when $\varepsilon^6 \leq L^2\sigma/n$.

## 1.6 Our Techniques

Let us first recall the main idea behind stochastic variance-reduced methods, such as SVRG [20].

The SVRG method divides iterations into epochs, each of length $n$. It maintains a snapshot point $\widetilde{\mathbf{x}}$ for each epoch, and computes the full gradient $\nabla f(\widetilde{\mathbf{x}})$ only for snapshots. Then, in each iteration $t$ at point $x_t$, SVRG defines gradient estimator $\widetilde{\nabla} = \nabla f_i(x_t) - \nabla f_i(\widetilde{\mathbf{x}}) + \nabla f(\widetilde{\mathbf{x}})$ which satisfies $\mathbb{E}_i[\widetilde{\nabla}] = \nabla f(x_t)$, and performs proximal update $x_{t+1} \leftarrow \mathsf{Prox}_{\psi,\alpha}(x_t - \alpha\widetilde{\nabla})$ for some learning rate $\alpha$. (Recall that if $\psi(\cdot) \equiv 0$ then we would have $x_{t+1} \leftarrow x_t - \alpha\widetilde{\nabla}$.)

In nearly all the aforementioned results for nonconvex optimization, researchers have either directly applied SVRG [6] (for the case $\sigma = L$), or repeatedly applied SVRG [1, 13] (for general $\sigma \in [0, L]$). This puts some limitation in the algorithmic design, because SVRG requires each epoch to be of length exactly $n$.[8]

**Our New Idea.** In this paper, we propose Natasha1 and Natasha1$^{\mathsf{full}}$ that are no longer black-box reductions to SVRG. Both of them still divide iterations into epochs of length $n$, and compute gradient estimators $\widetilde{\nabla}$ the same way as SVRG. However, we do not apply compute $x_t - \alpha\widetilde{\nabla}$ directly.

- In our base algorithm Natasha1, we divide each epoch into $p$ sub-epochs, each with a starting vector $\widehat{\mathbf{x}}$. Our theory suggests the choice $p \approx (\frac{\sigma^2}{L^2}n)^{1/3}$. Then, we replace the use of $\widetilde{\nabla}$ with $\widetilde{\nabla} + 2\sigma(x_t - \widehat{\mathbf{x}})$. This is equivalent to replacing $f(x)$ with its regularized version $f(x) + \sigma\|x - \widehat{\mathbf{x}}\|^2$, where the center $\widehat{\mathbf{x}}$ varies across sub-epochs. We provide pseudocode in Algorithm 1 and illustrate it in Figure 2.

  We view this additional term $2\sigma(x_t - \widehat{\mathbf{x}})$ as a type of **retraction**, which stabilizes the algorithm by moving the vector a bit in the backward direction towards $\widehat{\mathbf{x}}$.

- In our full algorithm Natasha1$^{\mathsf{full}}$, we add one more ingredient on top of Natasha1. That is, we perform updates $z_{t+1} \leftarrow \mathsf{Prox}_{\psi,\alpha}(z_t - \alpha\widetilde{\nabla})$ with respect to a different sequence $\{z_t\}$, and then define $x_t = \frac{1}{2}z_t + \frac{1}{2}\widehat{\mathbf{x}}$ and compute gradient estimators $\widetilde{\nabla}$ at points $x_t$. We provide pseudocode in Algorithm 2.

  We view this averaging $x_t = \frac{1}{2}z_t + \frac{1}{2}\widehat{\mathbf{x}}$ as another type of **retraction**, which stabilizes the algorithm by moving towards $\widehat{\mathbf{x}}$. The technique of computing gradients at points $x_t$ but moving a different sequence of points $z_t$ is related to the *Katyusha momentum* recently developed for convex optimization [2].

## 2 Preliminaries

Throughout this paper, we denote by $\|\cdot\|$ the Euclidean norm. We use $i \in_R [n]$ to denote that $i$ is generated from $[n] = \{1, 2, \ldots, n\}$ uniformly at random. We denote by $\nabla f(x)$ the full gradient of function $f$ if it is differentiable, and $\partial f(x)$ any subgradient if $f$ is only Lipschitz continuous at point $x$. We let $x^*$ be any minimizer of $F(x)$.

Recall some definitions on strong convexity (SC), nonconvexity, and smoothness.

---

[8]The epoch length of SVRG is always $n$ (or a constant multiple of $n$ in practice), because this ensures the computation of $\widetilde{\nabla}$ is of amortized gradient complexity $O(1)$. The per-iteration complexity of SVRG is thus the same as the traditional stochastic gradient descent (SGD).

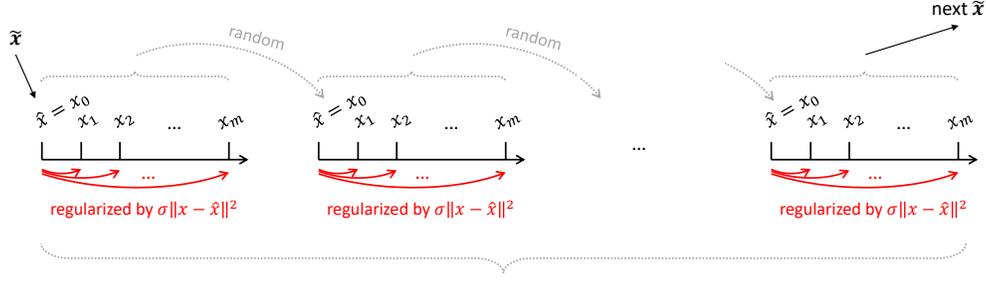

Figure 2: One full epoch of `Natasha1`. The $n$ iterations are divided into $p$ sub-epochs, each consisting of $m = n/p$ steps.

**Definition 2.1.** *For a function* $f \colon \mathbb{R}^d \to \mathbb{R}$,
- $f$ *is* σ-**strongly convex** *if* $\forall x, y \in \mathbb{R}^d$, *it satisfies* $f(y) \geq f(x) + \langle \partial f(x), y - x \rangle + \frac{\sigma}{2}\|x - y\|^2$.
- $f$ *is of* σ-*bounded nonconvexity (or* σ-**nonconvex** *for short) if* $\forall x, y \in \mathbb{R}^d$, *it satisfies*
$$f(y) \geq f(x) + \langle \partial f(x), y - x \rangle - \frac{\sigma}{2}\|x - y\|^2 \ .$$
- $f$ *is* $(\ell_1, \ell_2)$-**smooth** *if* $\forall x, y \in \mathbb{R}^d$, *it satisfies*
$$f(x) + \langle \nabla f(x), y - x \rangle + \tfrac{\ell_1}{2}\|x-y\|^2 \geq f(y) \geq f(x) + \langle \nabla f(x), y - x \rangle - \tfrac{\ell_2}{2}\|x-y\|^2 \ .$$
- $f$ *is* $L$-**smooth** *if it is* $(L, L)$-*smooth*.

The $(\ell_1, \ell_2)$-smoothness parameters were studied for instance by [12] to tackle the convex setting of problem (1.1).

**Definition 2.2.** *Given a parameter* $\eta > 0$, *the* **gradient mapping** *of* $F(\cdot)$ *in (1.1) at point* $x$ *is*
$$\mathcal{G}_\eta(x) \stackrel{\text{def}}{=} \frac{1}{\eta}(x - x') \qquad \text{where} \qquad x' = \arg\min_y \Big\{\psi(y) + \langle \nabla f(x), y \rangle + \frac{1}{2\eta}\|y - x\|^2\Big\}$$

*In particular, if* $\psi(\cdot) \equiv 0$, *then* $\mathcal{G}_\eta(x) \equiv \nabla f(x)$.

The following theorem for the SVRG method can be found in [12], which is built on the APPA/Catalyst reduction of [16, 24, 30]:

**Theorem 2.3** (SVRG). *Let* $G(y) \stackrel{\text{def}}{=} \psi(y) + \frac{1}{n}\sum_{i=1}^n g_i(y)$. *Suppose* $\psi(x)$ *is proper convex and* $g(x) = \frac{1}{n}\sum_{i=1}^n g_i(y)$ *is* σ-*strongly convex, then the SVRG method finds a point* $y$ *satisfying* $G(y) - G(y^*) \leq \varepsilon$

- *with gradient complexity* $O\big((n + \frac{L^2}{\sigma^2})\log\frac{1}{\varepsilon}\big)$, *if each* $g_i(\cdot)$ *is* $L$-*smooth (for* $L \geq \sigma$); *or*
- *with gradient complexity* $O\big((n + \frac{\ell_1 \ell_2}{\sigma^2})\log\frac{1}{\varepsilon}\big)$, *if each* $g_i(\cdot)$ *is* $(\ell_1, \ell_2)$-*smooth (for* $\ell_1, \ell_2 \geq \sigma$).

*If one performs the APPA/Catalyst reduction [16, 24, 30] to accelerate SVRG, the time gradient complexities improve to* $\widetilde{O}\big(n + n^{3/4}\sqrt{L/\sigma}\big)$ *and* $\widetilde{O}\big(n + \frac{n^{1/2}\sqrt{\ell_1+\ell_2}+n^{3/4}(\ell_1\ell_2)^{1/4}}{\sqrt{\sigma}}\big)$ *respectively.*[9]

After this paper appears online, a direct acceleration of SVRG (without using APPA/Catalyst) was obtained in [4]. In V1-V4 of this paper, we forgot an additive term $\frac{\sqrt{n(\ell_1+\ell_2)}}{\sqrt{\sigma}}$ in the complexity above. We have fixed this mistake in V5.

---

[9]Let us summarize the key idea of APPA/Catalyst for those readers not familiar with how it works. Suppose we are given a parameter $\sigma' \geq \sigma$ and an oracle that can, somehow, minimize $G(y) + \frac{\sigma'}{2}\|y - y_0\|^2$ for any arbitrary point $y_0$ in gradient complexity $T'$. Then, the result of APPA/Catalyst says that, we can use this oracle to minimize the original $G(y)$ in gradient complexity $T' \times \widetilde{O}(\sqrt{\sigma'/\sigma})$. In our case, let this oracle be SVRG. If each $g_i(x)$ is $L$-smooth, then we can choose $\sigma' = \sigma + L/\sqrt{n}$ and this ensures $T' = \widetilde{O}(n)$. Similarly, if each $g_i(x)$ is $(\ell_1, \ell_2)$-smooth, we can choose $\sigma' = \sigma + \sqrt{\ell_1 \ell_2/n} + (\ell_1 + \ell_2)/n$ and this ensures $T' = \widetilde{O}(n)$ as well.

## 2.1 RepeatSVRG

We recall the idea behind a simple algorithm [1, 13] —that we call repeatSVRG— which finds the $\varepsilon$-approximate stationary points for problem (1.1) when $f(x)$ is of $\sigma$-bounded nonconvexity. The algorithm is divided into stages. In each stage $t$, consider a modified function $F_t(x) \stackrel{\text{def}}{=} F(x) + \sigma\|x - x_t\|^2$. It is easy to see that $F_t(x)$ is $\sigma$-strongly convex, so one can apply the accelerated SVRG method to minimize $F_t(x)$. Let $x_{t+1}$ be any sufficiently accurate approximate minimizer of $F_t(x)$.[10]

Now, one can prove (see e.g. Section 4) that $x_{t+1}$ is an $O(\sigma\|x_t - x_{t+1}\|)$-approximate stationary point for $F(x)$. Therefore, if $\sigma\|x_t - x_{t+1}\| \leq \varepsilon$ we can stop the algorithm because we have already found an $O(\varepsilon)$-approximate stationary point. If $\sigma\|x_t - x_{t+1}\| > \varepsilon$, then it must satisfy that $F(x_t) - F(x_{t+1}) \geq \sigma\|x_t - x_{t+1}\|^2 \geq \Omega(\varepsilon^2/\sigma)$, but this cannot happen for more than $T = O\big(\frac{\sigma}{\varepsilon^2}(F(x_0) - F^*)\big)$ stages. Therefore, the total gradient complexity is $T$ multiplied with the complexity of accelerated SVRG in each stage, which is $\widetilde{O}(n + n^{3/4}\sqrt{L/\sigma})$ or $\widetilde{O}\big(n + \frac{n^{1/2}\sqrt{\ell_1+\ell_2}+n^{3/4}(\ell_1\ell_2)^{1/4}}{\sqrt{\sigma}}\big)$ according to Theorem 2.3.

*Remark* 2.4. The complexity of repeatSVRG can be inferred from [1, 13], but is not explicitly stated. For instance, the paper [13] repeatedly applies accelerated gradient descent (AGD) instead of accelerated SVRG; while the paper [1] repeatedly applies SVRG but did not state it as a separate lemma. Neither paper allows $F(x)$ to have a non-smooth proximal term $\psi(x)$. In any case, the credit of repeatSVRG belongs to the original authors of [1, 13], and our job here is only to state it properly in the most general setting.

In our revision of this paper, we became aware of two new results RapGrad [21] and Stagewise-Katyusha [14], which are different algorithms but in the similar spirit as repeatSVRG. The heart of these algorithms are also to modify the nonconvex objective $F(x)$ by adding a regularizer $\sigma\|x - x_t\|^2$, and then recursively solving the convex subproblem (but not using SVRG anymore). We acknowledge that RapGrad [21] is an independent work of repeatSVRG [1, 13].

## 3 Our Algorithms

We introduce two variants of our algorithms: (1) the base method `Natasha1` targets on the simple regime when $f(x)$ and each $f_i(x)$ are both $L$-smooth, and (2) the full method `Natasha1`[full] targets on the more refined regime when $f(x)$ is $L$-smooth but each $f_i(x)$ is $(\ell_1, \ell_2)$-smooth.

Both methods follow the general idea of variance-reduced stochastic gradient descent: in each inner-most iteration, they compute a *gradient estimator* $\widetilde{\nabla}$ that is of the form $\widetilde{\nabla} = \nabla f(\widetilde{\mathbf{x}}) - \nabla f_i(\widetilde{\mathbf{x}}) + \nabla f_i(x)$ and satisfies $\mathbb{E}_{i \in R[n]}[\widetilde{\nabla}] = \nabla f(x)$. Here, $\widetilde{\mathbf{x}}$ is a *snapshot* point that is changed once every $n$ iterations (i.e., for different $k = 1, 2, \ldots, T'$ in the pseudocode), and we call it a *full epoch* for every distinct $k$. Notice that the amortized gradient complexity for computing $\widetilde{\nabla}$ is $O(1)$ per-iteration.

**Base Method.** In `Natasha1` (see Algorithm 1), as illustrated by Figure 2, we divide each full epoch into $p$ sub-epochs $s = 0, 1, \ldots, p-1$, each of length $m = n/p$. In each sub-epoch $s$, we start with a point $x_0 = \widehat{x}$, and replace $f(x)$ with its regularized version $f^s(x) \stackrel{\text{def}}{=} f(x) + \sigma\|x - \widehat{x}\|^2$. Then, in each iteration $t$ of the sub-epoch $s$, we

- compute gradient estimator $\widetilde{\nabla}$ with respect to $f^s(x_t)$, and
- perform update $x_{t+1} = \arg\min_y \big\{\psi(y) + \langle\widetilde{\nabla}, y\rangle + \frac{1}{2\alpha}\|y - x_t\|^2\big\}$ with learning rate $\alpha$.

---
[10]Since the accelerated SVRG method has a linear convergence rate for strongly convex functions, the complexity to find such $x_{t+1}$ only depends logarithmically on this accuracy.

**Algorithm 1** `Natasha1`$(x^\varnothing, p, T', \alpha)$

**Input:** starting vector $x^\varnothing$, sub-epoch count $p \in [n]$, epoch count $T'$, learning rate $\alpha > 0$.
**Output:** vector $x^{\mathsf{out}}$.
1: $\widehat{\mathsf{x}} \leftarrow x^\varnothing$; $m \leftarrow n/p$; $X \leftarrow []$;
2: **for** $k \leftarrow 1$ **to** $T'$ **do**  ⋄ *$T'$ full epochs*
3:     $\widetilde{\mathsf{x}} \leftarrow \widehat{\mathsf{x}}$; $\mu \leftarrow \nabla f(\widetilde{\mathsf{x}})$;
4:     **for** $s \leftarrow 0$ **to** $p-1$ **do** ⋄ *$p$ sub-epochs in each epoch*
5:        $x_0 \leftarrow \widehat{\mathsf{x}}$; $X \leftarrow [X, \widehat{\mathsf{x}}]$;
6:        **for** $t \leftarrow 0$ **to** $m-1$ **do** ⋄ *$m$ iterations in each sub-epoch*
7:           $i \leftarrow$ a random choice from $\{1, \cdots, n\}$.
8:           $\widetilde{\nabla} \leftarrow \nabla f_i(x_t) - \nabla f_i(\widetilde{\mathsf{x}}) + \mu + 2\sigma(x_t - \widehat{\mathsf{x}})$ ⋄ $\mathbb{E}_i[\widetilde{\nabla}] = \nabla(f(x) + \sigma\|x - \widehat{\mathsf{x}}\|^2)\big|_{x_t}$
9:           $x_{t+1} = \arg\min_{y \in \mathbb{R}^d} \left\{ \psi(y) + \frac{1}{2\alpha}\|y - x_t\|^2 + \langle \widetilde{\nabla}, y \rangle \right\}$
10:        **end for**
11:        $\widehat{\mathsf{x}} \leftarrow$ a random choice from $\{x_0, x_1, \ldots, x_{m-1}\}$; ⋄ *for practitioners, choose the average*
12:     **end for**
13: **end for**
14: $\widehat{\mathsf{x}} \leftarrow$ a random vector in $X$; ⋄ *for practitioners, choose the last*
15: $x^{\mathsf{out}} \leftarrow$ an approximate minimizer of $G(y) \stackrel{\text{def}}{=} F(y) + \sigma\|y - \widehat{\mathsf{x}}\|^2$ using SVRG.
16: **return** $x^{\mathsf{out}}$. ⋄ *it suffices to run SVRG for $O(n \log \frac{L}{\varepsilon \sigma})$ iterations.*

---

**Algorithm 2** `Natasha1`$^{\mathsf{full}}(x^\varnothing, p, T', \alpha)$

**Input:** starting vector $x^\varnothing$, sub-epoch count $p \in [n]$, epoch count $T'$, learning rate $\alpha > 0$.
**Output:** vector $x^{\mathsf{out}}$.
1: $\widehat{\mathsf{x}} \leftarrow x^\varnothing$; $m \leftarrow n/p$; $X \leftarrow []$;
2: **for** $k \leftarrow 1$ **to** $T'$ **do** ⋄ *$T'$ full epochs*
3:     $\widetilde{\mathsf{x}} \leftarrow \widehat{\mathsf{x}}$; $\mu \leftarrow \nabla f(\widetilde{\mathsf{x}})$;
4:     **for** $s \leftarrow 0$ **to** $p-1$ **do** ⋄ *$p$ sub-epochs in each epoch*
5:        $z_0 \leftarrow \widehat{\mathsf{x}}$; $x_0 \leftarrow \widehat{\mathsf{x}}$; $X \leftarrow [X, \widehat{\mathsf{x}}]$;
6:        **for** $t \leftarrow 0$ **to** $m-1$ **do** ⋄ *$m$ iterations in each sub-epoch*
7:           $i \leftarrow$ a random choice from $\{1, \cdots, n\}$;
8:           $\widetilde{\nabla} \leftarrow \nabla f_i(x_t) - \nabla f_i(\widetilde{\mathsf{x}}) + \mu + 2\sigma(x_t - \widehat{\mathsf{x}})$; ⋄ $\mathbb{E}_i[\widetilde{\nabla}] = \nabla(f(x) + \sigma\|x - \widehat{\mathsf{x}}\|^2)\big|_{x_t}$
9:           $z_{t+1} = \arg\min_{y \in \mathbb{R}^d} \left\{ \psi(y) + \frac{1}{2\alpha}\|y - z_t\|^2 + \langle \widetilde{\nabla}, y \rangle \right\}$;
10:           $x_{t+1} = \frac{1}{2}z_{t+1} + \frac{1}{2}\widehat{\mathsf{x}}$; ⋄ *Katyusha momentum* $x_{t+1} = (1-\beta)z_{t+1} + \beta\widehat{\mathsf{x}}$
                                     ⋄ *theory predicts $\beta = \Theta\big(\frac{\sigma(\ell_1 + \ell_2)}{\ell_1 \ell_2}\big)$ gives the best performance*
                                     ⋄ *$\beta = 1/2$ however leads to the simplest proof*
11:        **end for**
12:        $\widehat{\mathsf{x}} \leftarrow$ a random choice from $\{x_0, x_1, \ldots, x_{m-1}\}$; ⋄ *for practitioners, choose the average*
13:     **end for**
14: **end for**
15: $\widehat{\mathsf{x}} \leftarrow$ a random vector in $X$; ⋄ *for practitioners, choose the last*
16: $x^{\mathsf{out}} \leftarrow$ an approximate minimizer of $G(y) \stackrel{\text{def}}{=} F(y) + \sigma\|y - \widehat{\mathsf{x}}\|^2$ using SVRG.
17: **return** $x^{\mathsf{out}}$. ⋄ *it suffices to run SVRG for $O(n \log \frac{L}{\varepsilon \sigma})$ iterations.*

Effectively, the introduction of the regularizer $\sigma\|x - \widehat{x}\|^2$ makes sure that when performing update $x_t \leftarrow x_{t+1}$, we also move a bit towards point $\widehat{x}$ (i.e., retraction by regularization). Finally, when the sub-epoch is done, we define $\widehat{x}$ to be a random one from $\{x_0, \ldots, x_{m-1}\}$.

**Full Method.** In `Natasha1`$^{\text{full}}$, we also divide each full epoch into $p$ sub-epochs. In each sub-epoch $s$, we start with a point $x_0 = z_0 = \widehat{x}$ and define $f^s(x) \stackrel{\text{def}}{=} f(x) + \sigma\|x - \widehat{x}\|^2$. However, this time in each iteration $t$, we

- compute gradient estimator $\widetilde{\nabla}$ with respect to $f^s(x_t)$,
- perform update $z_{t+1} = \arg\min_y \{\psi(y) + \langle \widetilde{\nabla}, y \rangle + \frac{1}{2\alpha}\|y - z_t\|^2\}$ with learning rate $\alpha$, and
- choose $x_{t+1} = \frac{1}{2}z_{t+1} + \frac{1}{2}\widehat{x}$.

Effectively, the regularizer $\sigma\|x - \widehat{x}\|^2$ makes sure that when performing updates, we move a bit towards point $\widehat{x}$ (i.e., retraction by regularization); at the same time, the choice $x_{t+1} = \frac{1}{2}z_{t+1} + \frac{1}{2}\widehat{x}$ also helps us move towards point $\widehat{x}$ (i.e., retraction by the so-called "Katyusha momentum"[11]). Finally, when the sub-epoch is over, we define $\widehat{x}$ to be a random one from the set $\{x_0, \ldots, x_{m-1}\}$, and move to the next sub-epoch.

## 4 A Sufficient Stopping Criterion

In this section, we present a sufficient condition for finding approximate stationary points in a $\sigma$-nonconvex function. Lemma 4.1 below states that, if we regularize the original function and define $G(x) \stackrel{\text{def}}{=} F(x) + \sigma\|x - \widehat{x}\|^2$ for an arbitrary point $\widehat{x}$, then the minimizer of $G(x)$ is an approximate saddle-point for $F(x)$.

**Lemma 4.1.** *Suppose $G(y) = F(y) + \sigma\|y - \widehat{x}\|^2$ for some given point $\widehat{x}$, and let $x^*$ be the unique minimizer of $G(y)$. If we minimize $G(y)$ and obtain a point $x$ satisfying*

$$G(x) - G(x^*) \leq \delta^2 \sigma ,$$

*then for every $\eta \in \left(0, \frac{1}{\max\{L, 4\sigma\}}\right]$, the gradient mapping $\mathcal{G}_\eta(x)$ of $F(x)$ satisfies*

$$\|\mathcal{G}_\eta(x)\|^2 \leq 12\sigma^2 \|x^* - \widehat{x}\|^2 + O\left(\frac{\delta^2}{\eta^2}\right) .$$

Notice that when $\psi(x) \equiv 0$ and $\delta = 0$, this lemma is trivial.[12] The main technical difficulty arises in order to deal with $\psi(x) \neq 0$ and the inexactness $\delta > 0$.

*Proof of Lemma 4.1.* Define auxiliary functions:

$$\Phi(y) \stackrel{\text{def}}{=} \psi(y) + \frac{1}{2\eta}\|y - x\|^2 + \langle \nabla f(x), y - x \rangle - \psi(x) \quad \text{and} \quad \overline{\Phi}(y) \stackrel{\text{def}}{=} \Phi(y) + \sigma\|y - \widehat{x}\|^2 - \sigma\|x - \widehat{x}\|^2$$

and let $z = \arg\min_y \Phi(y)$ and $\overline{z} = \arg\min_y \overline{\Phi}(y)$. In other words,

- $z$ is $x$ after applying a proximal update on $F(\cdot)$ with learning rate $\eta$ (so $\mathcal{G}_\eta(x) = \frac{x-z}{\eta}$), and
- $\overline{z}$ is $x$ after applying a proximal update on $G(\cdot)$ with learning rate $\eta$.

Observe that

---

[11]The idea for this second kind of retraction, and the idea of having the updates on a sequence $z_t$ but computing gradients at points $x_t$, is largely motivated by our recent work on the Katyusha momentum and the Katyusha acceleration [2].

[12]In this case, $\|\mathcal{G}_\eta(x)\| = \|\nabla F(x)\| = \|\nabla G(x) - 2\sigma(x - \widehat{x})\| = 2\sigma\|x - \widehat{x}\|$.

- $\Phi(\cdot)$ is $\frac{1}{\eta}$-strongly convex so $-\Phi(z) = \Phi(x) - \Phi(z) \geq \frac{1}{2\eta}\|z - x\|^2$;
- $\overline{\Phi}(\cdot)$ is $\frac{1}{\eta}$-strongly convex so $\overline{\Phi}(z) \geq \overline{\Phi}(\overline{z}) + \frac{1}{2\eta}\|z - \overline{z}\|^2$;
- $\overline{\Phi}(x) = 0 \geq \overline{\Phi}(\overline{z}) \geq G(\overline{z}) - G(x) \geq G(x^*) - G(x) \geq -\delta^2\sigma$ (since $\eta \leq 1/L$ and $f(\cdot)$ is $L$-smooth).

Summing the three inequalities up we have
$$\sigma\|z - \widehat{\mathsf{x}}\|^2 - \sigma\|x - \widehat{\mathsf{x}}\|^2 \geq -\delta^2\sigma + \frac{1}{2\eta}\|z - x\|^2 + \frac{1}{2\eta}\|z - \overline{z}\|^2 \ .$$

Since we have inequality $\|z - \widehat{\mathsf{x}}\|^2 = \|(z - \overline{z}) + (\overline{z} - \widehat{\mathsf{x}})\|^2 \leq (1 + 1/\beta)\|(z - \overline{z})\|^2 + (1 + \beta)\|(\overline{z} - \widehat{\mathsf{x}})\|^2$ for any $\beta > 0$, we can choose $\beta = 4\eta\sigma$ and obtain

$$(\sigma + \frac{1}{4\eta})\|(z - \overline{z})\|^2 + (\sigma + 4\eta\sigma^2)\|(\overline{z} - \widehat{\mathsf{x}})\|^2 - \sigma\|x - \widehat{\mathsf{x}}\|^2 \geq -\delta^2\sigma + \frac{1}{2\eta}\|z - x\|^2 + \frac{1}{2\eta}\|z - \overline{z}\|^2$$
$$\implies (\sigma + 4\eta\sigma^2)\|\overline{z} - \widehat{\mathsf{x}}\|^2 - \sigma\|x - \widehat{\mathsf{x}}\|^2 \geq -\delta^2\sigma + \frac{1}{2\eta}\|z - x\|^2 \quad (4.1)$$

where the implication uses the fact that $\frac{1}{4\eta} \geq \sigma$. At this point, notice that:

- We have $\|x - x^*\|^2 \leq \frac{2}{\sigma}(G(x) - G(x^*)) \leq 2\delta^2$ by the strong convexity of $G(\cdot)$, and thus using Young's inequality $-(1 + \beta)\|b\|^2 \leq -\|a + b\|^2 + (1 + 1/\beta)\|a\|^2$ for $\beta = \eta\sigma$, we have
$$-\sigma\|x - \widehat{\mathsf{x}}\|^2 \leq -(\sigma - \eta\sigma^2)\|x^* - \widehat{\mathsf{x}}\|^2 + O(\delta^2/\eta) \ .$$

- We have $\|\overline{z} - x^*\|^2 \leq \frac{2}{\sigma}(G(\overline{z}) - G(x^*)) \leq 2\delta^2$ because $G(\overline{z}) \leq G(x)$, and thus using Young's inequality $\|a + b\|^2 \leq (1 + \beta)\|b\|^2 + (1 + 1/\beta)\|a\|^2$ for $\beta = \eta\sigma$, we have
$$(\sigma + 4\eta\sigma^2)\|\overline{z} - \widehat{\mathsf{x}}\|^2 \leq (\sigma + 5\eta\sigma^2)\|x^* - \widehat{\mathsf{x}}\|^2 + O(\delta^2/\eta) \ .$$

Plugging them into (4.1), we have
$$(\sigma + 5\eta\sigma^2)\|x^* - \widehat{\mathsf{x}}\|^2 - (\sigma - \eta\sigma^2)\|x^* - \widehat{\mathsf{x}}\|^2 \geq \frac{1}{2\eta}\|z - x\|^2 - O(\delta^2/\eta) \ .$$

Rearranging it and multiplying both sides by $2/\eta$, we have
$$\|\mathcal{G}_\eta(x)\|^2 = \frac{1}{\eta^2}\|x - z\|^2 \leq 12\sigma^2\|x^* - \widehat{\mathsf{x}}\|^2 + O\Big(\frac{\delta^2}{\eta^2}\Big) \ . \qquad \square$$

## 5 Base Method: Analysis for One Full Epoch

In this section, we consider problem (1.1) where each $f_i(x)$ is $L$-smooth and $f(x)$ is $\sigma$-nonconvex. We use our base method `Natasha1` to minimize $F(x)$, and analyze its behavior for one full epoch in this section. We assume $\sigma \leq L$ without loss of generality, because any $L$-smooth function is also $L$-nonconvex.

**Notations.** We introduce the following notations for analysis purpose only.

- Let $\widehat{\mathsf{x}}^s$ be the vector $\widehat{\mathsf{x}}$ at the beginning of sub-epoch $s$.
- Let $x_t^s$ be the vector $x_t$ in sub-epoch $s$.
- Let $i_t^s$ be the index $i \in [n]$ in sub-epoch $s$ at iteration $t$.
- Let $f^s(x) \stackrel{\text{def}}{=} f(x) + \sigma\|x - \widehat{\mathsf{x}}^s\|^2$, $F^s(x) \stackrel{\text{def}}{=} F(x) + \sigma\|x - \widehat{\mathsf{x}}^s\|^2$, and $x_*^s \stackrel{\text{def}}{=} \arg\min_x\{F^s(x)\}$.
- Let $\widetilde{\nabla} f^s(x_t^s) \stackrel{\text{def}}{=} \nabla f_i(x_t^s) - \nabla f_i(\widetilde{\mathsf{x}}) + \nabla f(\widetilde{\mathsf{x}}) + 2\sigma(x_t - \widehat{\mathsf{x}})$ where $i = i_t^s$.

- Let $\widetilde{\nabla} f(x_t^s) \stackrel{\text{def}}{=} \nabla f_i(x_t^s) - \nabla f_i(\widetilde{\mathbf{x}}) + \nabla f(\widetilde{\mathbf{x}})$ where $i = i_t^s$.

We obviously have that $f^s(x)$ and $F^s(x)$ are $\sigma$-strongly convex, and $f^s(x)$ is $(L+2\sigma)$-smooth.

## 5.1 Variance Upper Bound

The following lemma gives an upper bound on the variance of the gradient estimator $\widetilde{\nabla} f^s(x_t^s)$:

**Lemma 5.1.** *We have* $\mathbb{E}_{i_t^s}\big[\|\widetilde{\nabla} f^s(x_t^s) - \nabla f^s(x_t^s)\|^2\big] \leq pL^2 \|x_t^s - \widehat{\mathbf{x}}^s\|^2 + pL^2 \sum_{k=0}^{s-1} \|\widehat{\mathbf{x}}^k - \widehat{\mathbf{x}}^{k+1}\|^2$ .

*Proof.* We have

$$\mathbb{E}_{i_t^s}\big[\|\widetilde{\nabla} f^s(x_t^s) - \nabla f^s(x_t^s)\|^2\big] = \mathbb{E}_{i_t^s}\big[\|\widetilde{\nabla} f(x_t^s) - \nabla f(x_t^s)\|^2\big]$$
$$= \mathbb{E}_{i \in_R [n]}\big[\|\big(\nabla f_i(x_t^s) - \nabla f_i(\widetilde{\mathbf{x}})\big) - \big(\nabla f(x_t^s) - \nabla f(\widetilde{\mathbf{x}})\big)\|^2\big]$$
$$\stackrel{\text{\textcircled{1}}}{\leq} \mathbb{E}_{i \in_R [n]}\big[\|\nabla f_i(x_t^s) - \nabla f_i(\widetilde{\mathbf{x}})\|^2\big]$$
$$\stackrel{\text{\textcircled{2}}}{\leq} p\mathbb{E}_{i \in_R [n]}\big[\|\nabla f_i(x_t^s) - \nabla f_i(\widehat{\mathbf{x}}^s)\|^2\big] + p\sum_{k=0}^{s-1} \mathbb{E}_{i \in_R [n]}\big[\|\nabla f_i(\widehat{\mathbf{x}}^k) - \nabla f_i(\widehat{\mathbf{x}}^{k+1})\|^2\big]$$
$$\stackrel{\text{\textcircled{3}}}{\leq} pL^2 \|x_t^s - \widehat{\mathbf{x}}^s\|^2 + pL^2 \sum_{k=0}^{s-1} \|\widehat{\mathbf{x}}^k - \widehat{\mathbf{x}}^{k+1}\|^2 .$$

Above, inequality ① is because for any random vector $\zeta \in \mathbb{R}^d$, it holds that $\mathbb{E}\|\zeta - \mathbb{E}\zeta\|^2 = \mathbb{E}\|\zeta\|^2 - \|\mathbb{E}\zeta\|^2$; inequality ② is because $\widehat{\mathbf{x}}^0 = \widetilde{\mathbf{x}}$ and for any $p$ vectors $a_1, a_2, \ldots, a_p \in \mathbb{R}^d$, it holds that $\|a_1 + \cdots + a_p\|^2 \leq p\|a_1\|^2 + \cdots + p\|a_p\|^2$; and inequality ③ is because each $f_i(\cdot)$ is $L$-smooth. □

## 5.2 Analysis for One Sub-Epoch

The following inequality is classically known as the "regret inequality" for mirror descent [11], and its proof is classical:

**Fact 5.2.** $\langle \widetilde{\nabla} f^s(x_t^s), x_{t+1}^s - u \rangle + \psi(x_{t+1}^s) - \psi(u) \leq \frac{\|x_t^s - u\|^2}{2\alpha} - \frac{\|x_{t+1}^s - u\|^2}{2\alpha} - \frac{\|x_{t+1}^s - x_t^s\|^2}{2\alpha}$ *for every* $u \in \mathbb{R}^d$.

*Proof.* Recall that the minimality of $x_{t+1}^s = \arg\min_{y \in \mathbb{R}^d}\{\frac{1}{2\alpha}\|y - x_t^s\|^2 + \psi(y) + \langle \widetilde{\nabla} f^s(x_t^s), y \rangle\}$ implies the existence of some subgradient $g \in \partial \psi(x_{t+1}^s)$ which satisfies $\frac{1}{\alpha}(x_{t+1}^s - x_t^s) + \widetilde{\nabla} f^s(x_t^s) + g = 0$. Combining this with $\psi(u) - \psi(x_{t+1}^s) \geq \langle g, u - x_{t+1}^s \rangle$, which is due to the convexity of $\psi(\cdot)$, we immediately have $\psi(u) - \psi(x_{t+1}^s) + \langle \frac{1}{\alpha}(x_{t+1}^s - x_t^s) + \widetilde{\nabla} f^s(x_t^s), u - x_{t+1}^s \rangle \geq \langle \frac{1}{\alpha}(x_{t+1}^s - x_t^s) + \widetilde{\nabla} f^s(x_t^s) + g, u - x_{t+1}^s \rangle = 0$. Rearranging this inequality we have

$$\langle \widetilde{\nabla} f^s(x_t^s), x_{t+1}^s - u \rangle + \psi(x_{t+1}^s) - \psi(u) \leq \langle -\frac{1}{\alpha}(x_{t+1}^s - x_t^s), x_{t+1}^s - u \rangle$$
$$= \frac{\|x_t^s - u\|^2}{2\alpha} - \frac{\|x_{t+1}^s - u\|^2}{2\alpha} - \frac{\|x_{t+1}^s - x_t^s\|^2}{2\alpha} . \quad \square$$

The following lemma is our main contribution for the base method `Natasha1`.

**Lemma 5.3.** *As long as* $\alpha \leq \frac{1}{2L+4\sigma}$, *we have*

$$\mathbb{E}\big[\big(F^s(\widehat{\mathbf{x}}^{s+1}) - F^s(x_*^s)\big)\big] \leq \mathbb{E}\Big[\frac{F^s(\widehat{\mathbf{x}}^s) - F^s(x_*^s)}{\sigma \alpha m/2} + \alpha pL^2\Big(\sum_{k=0}^{s} \|\widehat{\mathbf{x}}^k - \widehat{\mathbf{x}}^{k+1}\|^2\Big)\Big] .$$

*Proof.* We first compute that

$$F^s(x_{t+1}^s) - F^s(u) = f^s(x_{t+1}^s) - f^s(u) + \psi(x_{t+1}^s) - \psi(u)$$

$$\overset{①}{\leq} f^s(x_t^s) + \langle \nabla f^s(x_t^s), x_{t+1}^s - x_t^s \rangle + \frac{L+2\sigma}{2}\|x_t^s - x_{t+1}^s\|^2 - f^s(u) + \psi(x_{t+1}^s) - \psi(u)$$

$$\overset{②}{\leq} \langle \nabla f^s(x_t^s), x_{t+1}^s - x_t^s \rangle + \frac{L+2\sigma}{2}\|x_t^s - x_{t+1}^s\|^2 + \langle \nabla f^s(x_t^s), x_t^s - u \rangle + \psi(x_{t+1}^s) - \psi(u) \ . \quad (5.1)$$

Above, inequality ① uses the fact that $f^s(\cdot)$ is $(L+2\sigma)$-smooth; and inequality ② uses the convexity of $f^s(\cdot)$. Now, we take expectation with respect to $i_t^s$ on both sides of (5.1), and derive that:

$$\mathbb{E}_{i_t^s}\left[F^s(x_{t+1}^s)\right] - F^s(u)$$

$$\overset{①}{\leq} \mathbb{E}_{i_t^s}\left[\langle \widetilde{\nabla} f^s(x_t^s) - \nabla f^s(x_t^s), x_t^s - x_{t+1}^s \rangle + \langle \widetilde{\nabla} f^s(x_t^s), x_{t+1}^s - u \rangle + \frac{L+2\sigma}{2}\|x_t^s - x_{t+1}^s\|^2 + \psi(x_{t+1}^s) - \psi(u)\right]$$

$$\overset{②}{\leq} \mathbb{E}_{i_t^s}\left[\langle \widetilde{\nabla} f^s(x_t^s) - \nabla f^s(x_t^s), x_t^s - x_{t+1}^s \rangle + \frac{\|x_t^s - u\|^2}{2\alpha} - \frac{\|x_{t+1}^s - u\|^2}{2\alpha} - \left(\frac{1}{2\alpha} - \frac{L+2\sigma}{2}\right)\|x_{t+1}^s - x_t^s\|^2\right]$$

$$\overset{③}{\leq} \mathbb{E}_{i_t^s}\left[\alpha\|\widetilde{\nabla} f^s(x_t^s) - \nabla f^s(x_t^s)\|^2 + \frac{\|x_t^s - u\|^2}{2\alpha} - \frac{\|x_{t+1}^s - u\|^2}{2\alpha}\right]$$

$$\overset{④}{\leq} \mathbb{E}_{i_t^s}\left[\alpha p L^2\|x_t^s - \widehat{x}^s\|^2 + \alpha p L^2 \sum_{k=0}^{s-1}\|\widehat{x}^k - \widehat{x}^{k+1}\|^2 + \frac{\|x_t^s - u\|^2}{2\alpha} - \frac{\|x_{t+1}^s - u\|^2}{2\alpha}\right] \ . \quad (5.2)$$

Above, inequality ① is follows from (5.1) together with the fact that $\mathbb{E}_{i_t^s}[\widetilde{\nabla} f^s(x_t^s)] = \nabla f^s(x_t^s)$ implies

$$\mathbb{E}_{i_t^s}\left[\langle \nabla f^s(x_t^s), x_{t+1}^s - x_t^s \rangle + \langle \nabla f^s(x_t^s), x_t^s - u \rangle\right]$$
$$= \mathbb{E}_{i_t^s}\left[\langle \widetilde{\nabla} f^s(x_t^s) - \nabla f^s(x_t^s), x_t^s - x_{t+1}^s \rangle + \langle \widetilde{\nabla} f^s(x_t^s), x_{t+1}^s - u \rangle\right] \ ;$$

inequality ② uses Fact 5.2; inequality ③ uses $\alpha \leq \frac{1}{2L+4\sigma}$ together with Young's inequality $\langle a, b \rangle \leq \frac{1}{2}\|a\|^2 + \frac{1}{2}\|b\|^2$; and inequality ④ uses Lemma 5.1.

Finally, choosing $u = x_*^s$ to be the (unique) minimizer of $F^s(\cdot) = f^s(\cdot) + \psi(\cdot)$, and telescoping inequality (5.2) for $t = 0, 1, \ldots, m-1$, we have

$$\mathbb{E}\left[\sum_{t=1}^{m-1}\left(F^s(x_t^s) - F^s(x_*^s)\right)\right]$$

$$\leq \mathbb{E}\left[\frac{\|x_0^s - x_*^s\|^2}{2\alpha} + \sum_{t=0}^{m-1}\left(\alpha p L^2\|x_t^s - \widehat{x}^s\|^2 + \alpha p L^2 \sum_{k=0}^{s-1}\|\widehat{x}^k - \widehat{x}^{k+1}\|^2\right)\right]$$

$$\leq \mathbb{E}\left[\frac{F^s(\widehat{x}^s) - F^s(x_*^s)}{\sigma\alpha} + \alpha p m L^2\left(\sum_{k=0}^{s}\|\widehat{x}^k - \widehat{x}^{k+1}\|^2\right)\right] \ .$$

Above, the second inequality uses the fact that $\widehat{x}^{s+1}$ is chosen from $\{x_0^s, \ldots, x_{m-1}^s\}$ uniformly at random, as well as the $\sigma$-strong convexity of $F^s(\cdot)$.

Dividing both sides by $m$ and rearranging the terms (using $\frac{1}{2\sigma\alpha} \geq 1$), we have

$$\mathbb{E}\left[\left(F^s(\widehat{x}^{s+1}) - F^s(x_*^s)\right)\right] \leq \mathbb{E}\left[\frac{F^s(\widehat{x}^s) - F^s(x_*^s)}{\sigma\alpha m/2} + \alpha p L^2\left(\sum_{k=0}^{s}\|\widehat{x}^k - \widehat{x}^{k+1}\|^2\right)\right] \ . \quad \square$$

### 5.3 Analysis for One Full Epoch

One can telescope Lemma 5.3 for an entire epoch and arrive at the following lemma:

**Lemma 5.4.** *If $\alpha \leq \frac{1}{2L+4\sigma}$, $\alpha \geq \frac{4}{\sigma m}$ and $\alpha \leq \frac{\sigma}{p^2 L^2}$, we have*

$$\sum_{s=0}^{p-1} \mathbb{E}\Big[\big(F^s(\widehat{x}^s) - F^s(x_*^s)\big)\Big] \leq 2\mathbb{E}\Big[F(\widehat{x}^0) - F(\widehat{x}^p)\Big] \ .$$

*Proof.* Telescoping Lemma 5.3 for all the subepochs $s = 0, 1, \ldots, p-1$, we have

$$\sum_{s=0}^{p-1} \mathbb{E}\Big[\big(F^s(\widehat{x}^{s+1}) - F^s(x_*^s)\big)\Big] \leq \sum_{s=0}^{p-1} \mathbb{E}\Big[\frac{F^s(\widehat{x}^s) - F^s(x_*^s)}{\sigma \alpha m/2} + \alpha p^2 L^2 \|\widehat{x}^s - \widehat{x}^{s+1}\|^2\Big]$$

$$\overset{\text{\textcircled{1}}}{\leq} \sum_{s=0}^{p-1} \mathbb{E}\Big[\frac{F^s(\widehat{x}^s) - F^s(x_*^s)}{\sigma \alpha m/2} + \sigma \cdot \|\widehat{x}^{s+1} - \widehat{x}^s\|^2\Big]$$

$$\overset{\text{\textcircled{2}}}{=} \sum_{s=0}^{p-1} \mathbb{E}\Big[\frac{F^s(\widehat{x}^s) - F^s(x_*^s)}{\sigma \alpha m/2} + \big(F^s(\widehat{x}^{s+1}) - F^s(\widehat{x}^s)\big) - \big(F(\widehat{x}^{s+1}) - F(\widehat{x}^s)\big)\Big]$$

Above, ① uses $\alpha p^2 L^2 \leq \sigma$, and ② uses the definition $F^s(y) = F(y) + \sigma \|y - \widehat{x}^s\|^2$. Finally, rearranging both sides, and using the fact that $\frac{1}{\sigma \alpha m} \leq \frac{1}{4}$, we have the desired inequality. □

## 6 Base Method: Final Theorem

We are now ready to state and prove our main convergence theorem for `Natasha1`:

> **Theorem 1.** *Suppose in (1.1), each $f_i(x)$ is $L$-smooth and $f(x)$ is $\sigma$-nonconvex for $\sigma \leq L$. Then, if $\frac{L^2}{\sigma^2} \leq n$, $p = \Theta\big((\frac{\sigma^2}{L^2}n)^{1/3}\big)$ and $\alpha = \Theta(\frac{\sigma}{p^2 L^2})$, our base method `Natasha1` outputs a point $x^{\mathsf{out}}$ satisfying*
>
> $$\mathbb{E}[\|\mathcal{G}_\eta(x^{\mathsf{out}})\|^2] \leq O\Big(\frac{(L^2\sigma)^{1/3} n^{2/3}}{T'n}\Big) \cdot (F(x^\varnothing) - F^*) \ .$$
>
> *for $\eta = \frac{1}{\max\{L, 4\sigma\}}$. In other words, to obtain $\mathbb{E}[\|\mathcal{G}_\eta(x^{\mathsf{out}})\|^2] \leq \varepsilon^2$, we need gradient complexity*
>
> $$O\Big(n \log \frac{L}{\varepsilon \sigma} + \frac{(L^2\sigma)^{1/3} n^{2/3}}{\varepsilon^2} \cdot (F(x^\varnothing) - F^*)\Big) \ .$$

In the above theorem, we have assumed $\sigma \leq L$ without loss of generality because any $L$-smooth function is also $L$-nonconvex. Also, we have assumed $\frac{L^2}{\sigma^2} \leq n$ and if this inequality does not hold, then one should apply repeatSVRG for a faster running time (see Figure 1a).

*Proof of Theorem 1.* We choose $p = \big(\frac{\sigma^2}{24L^2}n\big)^{1/3}$, $m = n/p$, and $\alpha = \frac{4}{\sigma m} = \frac{\sigma}{6p^2 L^2} \leq \frac{1}{2L+4\sigma}$, so we can apply Lemma 5.4. If we telescope Lemma 5.4 for the entire algorithm (which has $T'$ full epochs), and use the fact that $\widehat{x}^p$ of the previous epoch equals $\widehat{x}^0$ of the next epoch, we conclude that if we choose a random epoch and a random subepoch $s$, we will have

$$\mathbb{E}[F^s(\widehat{x}^s) - F^s(x_*^s)] \leq \frac{2}{pT'}(F(x^\varnothing) - F^*) \ .$$

By the $\sigma$-strong convexity of $F^s(\cdot)$, we have $\mathbb{E}[\sigma \|\widehat{x}^s - x_*^s\|^2] \leq \frac{4}{pT'}(F(x^\varnothing) - F^*)$.

Now, $F^s(x) = F(x) + \sigma \|x - \widehat{x}^s\|^2$ satisfies the assumption of $G(x)$ in Lemma 4.1. If we use the SVRG method (see Theorem 2.3) to minimize the convex function $F^s(x)$, we get an output $x^{\mathsf{out}}$ satisfying $F^s(x^{\mathsf{out}}) - F^s(x_*^s) \leq \eta^2 \varepsilon^2 \sigma$ in gradient complexity $O\big((n + \frac{L^2}{\sigma^2}) \log \frac{1}{\eta \varepsilon \sigma}\big) \leq O(n \log \frac{L}{\varepsilon \sigma})$.

We can therefore apply Lemma 4.1 and conclude that this output $x^{\mathsf{out}}$ satisfies

$$\mathbb{E}[\|\mathcal{G}_\eta(x^{\mathsf{out}})\|^2] \leq O\Big(\frac{\sigma}{pT'}\Big) \cdot (F(x^\varnothing) - F^*) = O\Big(\frac{(L^2\sigma)^{1/3}n^{2/3}}{T'n}\Big) \cdot (F(x^\varnothing) - F^*) \ .$$

In other words, we obtain $\mathbb{E}[\|\mathcal{G}_\eta(x^{\mathsf{out}})\|^2] \leq \varepsilon^2$ using

$$T'n = O\Big(n + \frac{(L^2\sigma)^{1/3}n^{2/3}}{\varepsilon^2} \cdot (F(x^\varnothing) - F^*)\Big)$$

computations of the stochastic gradients. Here, the additive term $n$ is because $T' \geq 1$.

Finally, adding this with $O(n \log \frac{L}{\varepsilon\sigma})$, the gradient complexity for the application of SVRG in the last line of `Natasha1`, we finish the proof of the total gradient complexity. $\square$

## 7 Full Method: Analysis for One Full Epoch

In this section, we study a more refined version of problem (1.1), where $f(x)$ is $L$-smooth, each $f_i(x)$ is $(\ell_1, \ell_2)$-smooth, and $f(x)$ is $\sigma$-nonconvex. As later argued in Remark 8.1, we can assume $\sigma \leq \min\{\ell_1, \ell_2, L\}$ almost without loss of generality.

We use our full method `Natasha1`$^{\mathsf{full}}$ to minimize $F(x)$, and analyze its behavior for one full epoch in this section. Note that parameter $L$ is not needed in the specification of `Natasha1`$^{\mathsf{full}}$, but used only for analysis purpose.

**Notations.** We use the same notations as in Section 5, with an additional one highlighted here:

- Let $\widehat{x}^s$ be the vector $\widehat{x}$ at the beginning of sub-epoch $s$.
- Let $x_t^s$ be the vector $x_t$ in sub-epoch $s$.
- Let $i_t^s$ be the index $i \in [n]$ in sub-epoch $s$ at iteration $t$.
- Let $F^s(x) \stackrel{\text{def}}{=} F(x) + \sigma\|x - \widehat{x}^s\|^2$ and $x_*^s \stackrel{\text{def}}{=} \arg\min_x\{F^s(x)\}$.
- Let $f^s(x) \stackrel{\text{def}}{=} f(x) + \sigma\|x - \widehat{x}^s\|^2$ and $\boxed{f_i^s(x) \stackrel{\text{def}}{=} f_i(x) + \sigma\|x - \widehat{x}^s\|^2}$.
- Let $\widetilde{\nabla} f^s(x_t^s) \stackrel{\text{def}}{=} \nabla f_i(x_t^s) - \nabla f_i(\widetilde{x}) + \nabla f(\widetilde{x}) + 2\sigma(x_t - \widehat{x})$ where $i = i_t^s$.
- Let $\widetilde{\nabla} f(x_t^s) \stackrel{\text{def}}{=} \nabla f_i(x_t^s) - \nabla f_i(\widetilde{x}) + \nabla f(\widetilde{x})$ where $i = i_t^s$.

We obviously have that $f^s(x)$ and $F^s(x)$ are $\sigma$-strongly convex, and $f^s(x)$ is $(L + 2\sigma)$-smooth.

### 7.1 Variance Upper Bound

In this subsection we derive a new upper bound on the variance of the gradient estimator $\widetilde{\nabla}$. This bound will be tighter than Lemma 5.1, and will make use of the asymmetry between parameters $\ell_1$ and $\ell_2$. To achieve so, we first need to introduce the following lemma:

**Lemma 7.1.** *If $g(y) = \frac{1}{n}\sum_{i=1}^n g_i(y)$ is convex, and if each $g_i$ is $(\ell_1, \ell_2)$-smooth, then we have*

$$\mathbb{E}_{i \in_R [n]}\big[\|\nabla g_i(y_1) - \nabla g_i(y_2)\|^2\big]$$
$$\leq 2(\ell_1 + \ell_2)(g(y_2) - g(y_1) - \langle \nabla g(y_1), y_2 - y_1\rangle)] + 6\ell_1\ell_2\|y_2 - y_1\|^2 \ .$$

*Proof.* We consider two cases: $\ell_2 \leq \ell_1$ and $\ell_2 \geq \ell_1$.

- In the first case, we define $\phi_i(z) \stackrel{\text{def}}{=} g_i(z) - \langle \nabla g_i(y_1), z\rangle + \frac{\ell_2}{2}\|z - y_1\|^2$ for each $i \in [n]$. This function $\phi_i(z)$ is a convex, $(\ell_1 + \ell_2)$-smooth function that has a minimizer $z = y_1$ (which can

be seen by taking the derivative). For this reason, we claim that

$$\forall z: \quad \phi_i(y_1) \leq \phi_i(z) - \frac{1}{\ell_1 + \ell_2}\|\nabla \phi_i(z)\|^2 \ , \tag{7.1}$$

and this inequality is classical for smooth functions (see for instance Theorem 2.1.5 in textbook [25]). By expanding out the definition of $\phi_i(\cdot)$ in (7.1), we immediately have

$$g_i(y_1) - \langle \nabla g_i(y_1), y_1 \rangle \leq g_i(z) - \langle \nabla g_i(y_1), z \rangle + \frac{\ell_2}{2}\|z - y_1\|^2$$
$$- \frac{1}{2(\ell_1 + \ell_2)}\|\nabla g_i(z) - \nabla g_i(y_1) + \ell_2(z - y_1)\|^2$$

which then implies

$$\|\nabla g_i(z) - \nabla g_i(y_1)\|^2 \leq 2\|\nabla g_i(z) - \nabla g_i(y_1) + \ell_2(z - y_1)\|^2 + 2\|\ell_2(z - y_1)\|^2$$
$$\leq 2(\ell_1 + \ell_2)(g_i(z) - g_i(y_1) - \langle \nabla g_i(y_1), z - y_1 \rangle) + (4\ell_2^2 + 2\ell_1\ell_2)\|z - y_1\|^2 \ . \tag{7.2}$$

Now, by choosing $z = y_2$ and taking expectation with $i$ in (7.2), we have

$$\mathbb{E}_i\big[\|\nabla g_i(y_2) - \nabla g_i(y_1)\|^2\big]$$
$$\leq 2(\ell_1 + \ell_2)\big(g(y_2) - g(y_1) - \langle \nabla g(y_1), y_2 - y_1 \rangle\big) + (4\ell_2^2 + 2\ell_1\ell_2)\|y_2 - y_1\|^2 \tag{7.3}$$

- In the second case, we define $\phi_i(z) \stackrel{\text{def}}{=} -g_i(z) + \langle \nabla g_i(y_2), z \rangle + \frac{\ell_1}{2}\|z - y_2\|^2$ for each $i \in [n]$. It is clear that $\phi_i(z)$ is a convex, $(\ell_1 + \ell_2)$-smooth function that has a minimizer $z = y_2$ (which can be seen by taking the derivative). For this reason, we have

$$\forall z: \quad \phi_i(y_2) \leq \phi_i(z) - \frac{1}{\ell_1 + \ell_2}\|\nabla \phi_i(z)\|^2 \ . \tag{7.4}$$

By expanding out the definition of $\phi_i(\cdot)$ in (7.4), we immediately have

$$- g_i(y_2) + \langle \nabla g_i(y_2), y_2 \rangle \leq -g_i(z) + \langle \nabla g_i(y_2), z \rangle + \frac{\ell_1}{2}\|z - y_2\|^2$$
$$- \frac{1}{2(\ell_1 + \ell_2)}\|\nabla g_i(z) - \nabla g_i(y_2) - \ell_1(z - y_2)\|^2$$

which then implies that

$$\|\nabla g_i(z) - \nabla g_i(y_2)\|^2 \leq 2\|\nabla g_i(z) - \nabla g_i(y_2) - \ell_1(z - y_2)\|^2 + 2\|\ell_2(z - y_2)\|^2$$
$$\leq 2(\ell_1 + \ell_2)(g_i(y_2) - g_i(z) + \langle \nabla g_i(y_2), z - y_2 \rangle) + (4\ell_1^2 + 2\ell_1\ell_2)\|z - y_2\|^2 \ . \tag{7.5}$$

Now by choosing $z = y_1$ and taking expectation over $i$ in (7.5), we have

$$\mathbb{E}_i\big[\|\nabla g_i(y_1) - \nabla g_i(y_2)\|^2\big]$$
$$\leq 2(\ell_1 + \ell_2)\big(g(y_2) - g(y_1) + \langle \nabla g(y_2), y_1 - y_2 \rangle\big) + (4\ell_1^2 + 2\ell_1\ell_2)\|y_1 - y_2\|^2$$
$$\leq (4\ell_1^2 + 2\ell_1\ell_2)\|y_1 - y_2\|^2$$
$$\leq 2(\ell_1 + \ell_2)\big(g(y_2) - g(y_1) - \langle \nabla g(y_1), y_2 - y_1 \rangle\big) + (4\ell_2^2 + 2\ell_1\ell_2)\|y_2 - y_1\|^2 \ . \tag{7.6}$$

Above, the second and third inequalities use the convexity of $g(\cdot)$.

Combining (7.3) and (7.6) we finish the proof of the lemma. □

We are now ready to state our final variance upper bound:

**Lemma 7.2** (variance bound). *There exists constant $C \geq 1$ such that, if we define*

- $\Phi^s(y) \stackrel{\text{def}}{=} C(\ell_1 + \ell_2) \cdot (f^s(\widehat{x}^s) - f^s(y) - \langle \nabla f^s(y), \widehat{x}^s - y \rangle)] + C(\ell_1 \ell_2) \cdot \|y - \widehat{x}^s\|^2 \geq 0;$
- $\Phi_t^s = \Phi^s(x_t^s)$ and $\Phi^s = \Phi^s(\widehat{x}^{s+1})$,

*then, we have* $\mathbb{E}_i\big[\|\widetilde{\nabla} f^s(x_t^s) - \nabla f^s(x_t^s)\|^2\big] \leq p\Phi_t^s + p \sum_{k=0}^{s-1} \Phi^k$ *where* $i = i_t^s$.

Before proceeding to the proof, we point out that if $\ell_1 = \ell_2 = L$ like in the base setting, then we shall have $\Phi^s(y) \leq O(L^2)\|y - \widehat{x}^s\|^2$ and Lemma 7.2 becomes identical to Lemma 5.1.

*Proof.* If we plug in $g = f^s$ and $g_i = f_i^s$ in Lemma 7.1, we have $g_i$ is $(\ell_1 + 2\sigma, \ell_2 - 2\sigma)$-smooth and thus each $g_i$ is also $(3\ell_1, \ell_2)$-smooth. Therefore, Lemma 7.1 implies there exists constant $C \geq 1$ such that

$$\mathbb{E}_i\big[\|\nabla f_i(y) - \nabla f_i(\widehat{x}^s)\|^2\big] \leq 2\mathbb{E}_i\big[\|\nabla f_i^s(y) - \nabla f_i^s(\widehat{x}^s)\|^2\big] + 2\|2\sigma(y - \widehat{x}^s)\|^2$$
$$\leq C(\ell_1 + \ell_2) \cdot (f^s(\widehat{x}^s) - f^s(y) - \langle \nabla f^s(y), \widehat{x}^s - y \rangle)] + C(\ell_1 \ell_2) \cdot \|y - \widehat{x}^s\|^2$$
$$= \Phi^s(y) \ . \tag{7.7}$$

Therefore, the variance term:

$$\mathbb{E}_i\big[\|\widetilde{\nabla} f^s(x_t^s) - \nabla f^s(x_t^s)\|^2\big] = \mathbb{E}_i\big[\|\widetilde{\nabla} f(x_t^s) - \nabla f(x_t^s)\|^2\big]$$
$$= \mathbb{E}_i\big[\|\big(\nabla f_i(x_t^s) - \nabla f_i(\widetilde{x})\big) - \big(\nabla f(x_t^s) - \nabla f(\widetilde{x})\big)\|^2\big]$$
$$\stackrel{\text{①}}{\leq} \mathbb{E}_i\big[\|\nabla f_i(x_t^s) - \nabla f_i(\widetilde{x})\|^2\big]$$
$$\stackrel{\text{②}}{\leq} p\mathbb{E}_i\big[\|\nabla f_i(x_t^s) - \nabla f_i(\widehat{x}^s)\|^2\big] + p\sum_{k=0}^{s-1}\mathbb{E}_i\big[\|\nabla f_i(\widehat{x}^k) - \nabla f_i(\widehat{x}^{k+1})\|^2\big]$$
$$\stackrel{\text{③}}{\leq} p\Phi_t^s + p\sum_{k=0}^{s-1} \Phi^k \ . \tag{7.8}$$

Above, inequality ① is because for any random vector $\zeta \in \mathbb{R}^d$, it holds that $\mathbb{E}\|\zeta - \mathbb{E}\zeta\|^2 = \mathbb{E}\|\zeta\|^2 - \|\mathbb{E}\zeta\|^2$; inequality ② is because $\widehat{x}^0 = \widetilde{x}$ and for any $p$ vectors $a_1, a_2, \ldots, a_p \in \mathbb{R}^d$, it holds that $\|a_1 + \cdots + a_p\|^2 \leq p\|a_1\|^2 + \cdots + p\|a_p\|^2$; and inequality ③ is from repeatedly applying (7.7). □

## 7.2 Analysis for One Sub-Epoch

The following fact is analogous to Fact 5.2, and the only difference is that in `Natasha1`[full] we are applying proximal updates on the $\{z_t^s\}_t$ sequence.

**Fact 7.3.** $\langle \widetilde{\nabla} f^s(x_t^s), z_{t+1}^s - u \rangle + \psi(z_{t+1}^s) - \psi(u) \leq \frac{\|z_t^s - u\|^2}{2\alpha} - \frac{\|z_{t+1}^s - u\|^2}{2\alpha} - \frac{\|z_{t+1}^s - z_t^s\|^2}{2\alpha}$ *for every* $u \in \mathbb{R}^d$.

*Proof.* Recall that the minimality of $z_{t+1}^s = \arg\min_{y \in \mathbb{R}^d}\{\frac{1}{2\alpha}\|y - z_t^s\|^2 + \psi(y) + \langle \widetilde{\nabla} f^s(x_t^s), y \rangle\}$ implies the existence of some subgradient $g \in \partial \psi(z_{t+1}^s)$ which satisfies $\frac{1}{\alpha}(z_{t+1}^s - z_t^s) + \widetilde{\nabla} f^s(x_t^s) + g = 0$. Combining this with $\psi(u) - \psi(z_{t+1}^s) \geq \langle g, u - z_{t+1}^s \rangle$, which is due to the convexity of $\psi(\cdot)$, we immediately have $\psi(u) - \psi(z_{t+1}^s) + \langle \frac{1}{\alpha}(z_{t+1}^s - z_t^s) + \widetilde{\nabla} f^s(x_t^s), u - z_{t+1}^s \rangle \geq \langle \frac{1}{\alpha}(z_{t+1}^s - z_t^s) + \widetilde{\nabla} f^s(x_t^s) + g, u - z_{t+1}^s \rangle = 0$. Rearranging this inequality we have

$$\langle \widetilde{\nabla} f^s(x_t^s), z_{t+1}^s - u \rangle + \psi(z_{t+1}^s) - \psi(u) \leq \langle -\frac{1}{\alpha}(z_{t+1}^s - z_t^s), z_{t+1}^s - u \rangle$$
$$= \frac{\|z_t^s - u\|^2}{2\alpha} - \frac{\|z_{t+1}^s - u\|^2}{2\alpha} - \frac{\|z_{t+1}^s - z_t^s\|^2}{2\alpha} \ . \quad □$$

The following lemma is our technical main contribution for the full method `Natasha1`[full].

**Lemma 7.4.** *If $\alpha \leq \frac{1}{L+2\sigma}$, then we have the following inequality for sub-epoch $s$:*

$$\mathbb{E}\Big[\big(F^s(\widehat{x}^{s+1}) - F^s(x_*^s)\big)\Big]$$
$$\leq \mathbb{E}\Big[\frac{F^s(\widehat{x}^s) - F^s(x_*^s)}{\sigma\alpha m/2} + \alpha p\Big(\sum_{k=0}^{s}\Phi^k\Big) + \langle \nabla f^s(\widehat{x}^{s+1}), \widehat{x}^s - \widehat{x}^{s+1}\rangle + \big(\psi(\widehat{x}^s) - \psi(\widehat{x}^{s+1})\big)\Big] \ .$$

*Proof.* We first compute that

$$2F^s(x_{t+1}^s) - F^s(x_t^s) - F^s(u) = 2f^s(x_{t+1}^s) - f^s(x_t^s) - f^s(u) + 2\psi(x_{t+1}^s) - \psi(x_t^s) - \psi(u)$$
$$\overset{①}{\leq} f^s(x_t^s) + 2\langle \nabla f^s(x_t^s), x_{t+1}^s - x_t^s\rangle + (L+2\sigma)\|x_t^s - x_{t+1}^s\|^2 - f^s(u) + 2\psi(x_{t+1}^s) - \psi(x_t^s) - \psi(u)$$
$$\overset{②}{=} f^s(x_t^s) + \langle \nabla f^s(x_t^s), z_{t+1}^s - z_t^s\rangle + \frac{L+2\sigma}{4}\|z_t^s - z_{t+1}^s\|^2 - f^s(u) + 2\psi(x_{t+1}^s) - \psi(x_t^s) - \psi(u)$$
$$\overset{③}{\leq} \langle \nabla f^s(x_t^s), z_{t+1}^s - z_t^s\rangle + \frac{L+2\sigma}{4}\|z_t^s - z_{t+1}^s\|^2 + \langle \nabla f^s(x_t^s), x_t^s - u\rangle + \psi(z_{t+1}^s) + \psi(\widehat{x}^s) - \psi(x_t^s) - \psi(u)$$
(7.9)

Above, inequality ① uses the fact that $f^s(\cdot)$ is $(L+2\sigma)$-smooth; equality ② uses the fact that $z_{t+1}^s - z_t^s = 2(x_{t+1}^s - x_t^s)$; inequality ③ uses the convexity of $f^s(\cdot)$, the convexity of $\psi(\cdot)$, and the fact $x_{t+1}^s = \frac{1}{2}(z_{t+1}^s + \widehat{x}^s)$.

Now, we take expectation with respect to $i_t^s$ on both sides of (7.9), and derive that:

$$2\mathbb{E}_{i_t^s}\big[F^s(x_{t+1}^s)\big] - F^s(x_t^s) - F^s(u)$$
$$\overset{①}{\leq} \mathbb{E}_{i_t^s}\Big[\langle \widetilde{\nabla} f^s(x_t^s) - \nabla f^s(x_t^s), z_t^s - z_{t+1}^s\rangle + \langle \widetilde{\nabla} f^s(x_t^s), z_{t+1}^s - u\rangle + \frac{L+2\sigma}{4}\|z_t^s - z_{t+1}^s\|^2 + \psi(z_{t+1}^s) - \psi(u)\Big]$$
$$\qquad + \langle \nabla f^s(x_t^s), x_t^s - z_t^s\rangle + \psi(\widehat{x}^s) - \psi(x_t^s)$$
$$\overset{②}{\leq} \mathbb{E}_{i_t^s}\Big[\langle \widetilde{\nabla} f^s(x_t^s) - \nabla f^s(x_t^s), z_t^s - z_{t+1}^s\rangle + \frac{\|z_t^s - u\|^2}{2\alpha} - \frac{\|z_{t+1}^s - u\|^2}{2\alpha} - \Big(\frac{1}{2\alpha} - \frac{L+2\sigma}{4}\Big)\|z_{t+1}^s - z_t^s\|^2\Big]$$
$$\qquad + \langle \nabla f^s(x_t^s), x_t^s - z_t^s\rangle + \psi(\widehat{x}^s) - \psi(x_t^s)$$
$$\overset{③}{\leq} \mathbb{E}_{i_t^s}\Big[\alpha\|\widetilde{\nabla} f^s(x_t^s) - \nabla f^s(x_t^s)\|^2 + \frac{\|z_t^s - u\|^2}{2\alpha} - \frac{\|z_{t+1}^s - u\|^2}{2\alpha}\Big] + \langle \nabla f^s(x_t^s), x_t^s - z_t^s\rangle + \psi(\widehat{x}^s) - \psi(x_t^s)$$
$$\overset{④}{\leq} \mathbb{E}_{i_t^s}\Big[\alpha p\Phi_t^s + \alpha p\sum_{k=0}^{s-1}\Phi^k + \frac{\|z_t^s - u\|^2}{2\alpha} - \frac{\|z_{t+1}^s - u\|^2}{2\alpha}\Big] + \langle \nabla f^s(x_t^s), \widehat{x}^s - x_t^s\rangle + \psi(\widehat{x}^s) - \psi(x_t^s) \ .$$

Above, inequality ① is from (7.9) together with the fact that $\mathbb{E}_{i_t^s}[\widetilde{\nabla} f^s(x_t^s)] = \nabla f^s(x_t^s)$ implies

$$\mathbb{E}_{i_t^s}\big[\langle \nabla f^s(x_t^s), z_{t+1}^s - z_t^s\rangle + \langle \nabla f^s(x_t^s), x_t^s - u\rangle\big]$$
$$= \mathbb{E}_{i_t^s}\big[\langle \nabla f^s(x_t^s), x_t^s - z_t^s\rangle + \langle \widetilde{\nabla} f^s(x_t^s) - \nabla f^s(x_t^s), z_t^s - z_{t+1}^s\rangle + \langle \widetilde{\nabla} f^s(x_t^s), z_{t+1}^s - u\rangle\big] \ ;$$

inequality ② uses Fact 7.3; inequality ③ uses $\alpha \leq \frac{1}{L+2\sigma}$ together with Young's inequality $\langle a, b\rangle \leq \frac{1}{2}\|a\|^2 + \frac{1}{2}\|b\|^2$; and inequality ④ uses Lemma 7.2 and the fact that $x_t^s = \frac{1}{2}z_t^s + \frac{1}{2}\widehat{x}^s$.

Finally, choosing $u = x_*^s$ to be the (unique) minimizer of $F^s(\cdot) = f^s(\cdot) + \psi(\cdot)$, and telescoping the above inequality for $t = 0, 1, \ldots, m-1$, we have

$$\mathbb{E}\Big[\sum_{t=1}^{m-1}\big(F^s(x_t^s) - F^s(x_*^s)\big)\Big]$$
$$\leq \mathbb{E}\Big[\frac{\|z_0^s - x_*^s\|^2}{2\alpha} + \sum_{t=0}^{m-1}\Big(\alpha p\Phi_t^s + \alpha p\sum_{k=0}^{s-1}\Phi^k + \langle \nabla f^s(x_t^s), \widehat{x}^s - x_t^s\rangle + \psi(\widehat{x}^s) - \psi(x_t^s)\Big)\Big] \ .$$

Using the fact $\widehat{x}^{s+1}$ is chosen uniformly at random from $\{x_0^s, \ldots, x_{m-1}^s\}$, and the fact that $x_0^s = \widehat{x}^s$, the above inequality implies

$$\mathbb{E}\Big[m\big(F^s(\widehat{x}^{s+1}) - F^s(x_*^s)\big) - \big(F^s(\widehat{x}^s) - F^s(x_*^s)\big)\Big]$$

$$\leq \mathbb{E}\Big[\frac{\|z_0^s - x_*^s\|^2}{2\alpha} + \alpha pm\Big(\sum_{k=0}^{s}\Phi^k\Big) + m\langle\nabla f^s(\widehat{x}^{s+1}), \widehat{x}^s - \widehat{x}^{s+1}\rangle + m\big(\psi(\widehat{x}^s) - \psi(\widehat{x}^{s+1})\big)\Big]$$

$$\leq \mathbb{E}\Big[\frac{F^s(\widehat{x}^s) - F^s(x_*^s)}{\sigma\alpha} + \alpha pm\Big(\sum_{k=0}^{s}\Phi^k\Big) + m\langle\nabla f^s(\widehat{x}^{s+1}), \widehat{x}^s - \widehat{x}^{s+1}\rangle + m\big(\psi(\widehat{x}^s) - \psi(\widehat{x}^{s+1})\big)\Big] \ .$$

Above, the second inequality uses the fact that $z_0^s = \widehat{x}^s$ and that $F^s(\cdot)$ is $\sigma$-strongly convex. Dividing both sides by $m$ and rearranging the terms (using $\frac{1}{\sigma\alpha} \geq 1$), we have

$$\mathbb{E}\Big[\big(F^s(\widehat{x}^{s+1}) - F^s(x_*^s)\big)\Big]$$

$$\leq \mathbb{E}\Big[\frac{F^s(\widehat{x}^s) - F^s(x_*^s)}{\sigma\alpha m/2} + \alpha p\Big(\sum_{k=0}^{s}\Phi^k\Big) + \langle\nabla f^s(\widehat{x}^{s+1}), \widehat{x}^s - \widehat{x}^{s+1}\rangle + \big(\psi(\widehat{x}^s) - \psi(\widehat{x}^{s+1})\big)\Big] \ . \qquad \square$$

### 7.3 Analysis for One Full Epoch

We telescope Lemma 7.4 for an entire epoch and arrive at the following lemma:

**Lemma 7.5.** *If $\alpha \leq O(\frac{\sigma}{p^2\ell_1\ell_2})$ and $\alpha \geq \Omega(\frac{1}{\sigma m})$, we have*

$$\sum_{s=0}^{p-1}\mathbb{E}\Big[\big(F^s(\widehat{x}^s) - F^s(x_*^s)\big)\Big] \leq 3\mathbb{E}\Big[F(\widehat{x}^0) - F(\widehat{x}^p)\Big] \ .$$

*Proof.* Telescoping Lemma 7.4 for all the subepochs $s = 0, 1, \ldots, p-1$, we have

$$\sum_{s=0}^{p-1}\mathbb{E}\Big[\big(F^s(\widehat{x}^{s+1}) - F^s(x_*^s)\big)\Big]$$

$$\overset{①}{\leq} \sum_{s=0}^{p-1}\mathbb{E}\Big[\frac{F^s(\widehat{x}^s) - F^s(x_*^s)}{\sigma\alpha m/2} + \alpha p^2\Phi^s + \langle\nabla f^s(\widehat{x}^{s+1}), \widehat{x}^s - \widehat{x}^{s+1}\rangle + \big(\psi(\widehat{x}^s) - \psi(\widehat{x}^{s+1})\big)\Big]$$

$$\overset{②}{\leq} \sum_{s=0}^{p-1}\mathbb{E}\Big[\frac{F^s(\widehat{x}^s) - F^s(x_*^s)}{\sigma\alpha m/2} + \langle\nabla f^s(\widehat{x}^{s+1}), \widehat{x}^s - \widehat{x}^{s+1}\rangle + \big(\psi(\widehat{x}^s) - \psi(\widehat{x}^{s+1})\big)$$
$$+ \alpha p^2 C(\ell_1 + \ell_2)\cdot(f^s(\widehat{x}^s) - f^s(\widehat{x}^{s+1}) - \langle\nabla f^s(\widehat{x}^{s+1}), \widehat{x}^s - \widehat{x}^{s+1}\rangle)\Big] + \alpha p^2 C(\ell_1\ell_2)\cdot\|\widehat{x}^{s+1} - \widehat{x}^s\|^2\Big]$$

$$\overset{③}{\leq} \sum_{s=0}^{p-1}\mathbb{E}\Big[\frac{F^s(\widehat{x}^s) - F^s(x_*^s)}{\sigma\alpha m/2} + \langle\nabla f^s(\widehat{x}^{s+1}), \widehat{x}^s - \widehat{x}^{s+1}\rangle + \big(\psi(\widehat{x}^s) - \psi(\widehat{x}^{s+1})\big)$$
$$+ (f^s(\widehat{x}^s) - f^s(\widehat{x}^{s+1}) - \langle\nabla f^s(\widehat{x}^{s+1}), \widehat{x}^s - \widehat{x}^{s+1}\rangle)\Big] + 2\sigma\cdot\|\widehat{x}^{s+1} - \widehat{x}^s\|^2\Big]$$

$$= \sum_{s=0}^{p-1}\mathbb{E}\Big[\frac{F^s(\widehat{x}^s) - F^s(x_*^s)}{\sigma\alpha m/2} + \big(F^s(\widehat{x}^s) - F^s(\widehat{x}^{s+1})\big) + 2\sigma\cdot\|\widehat{x}^{s+1} - \widehat{x}^s\|^2\Big]$$

$$\overset{④}{=} \sum_{s=0}^{p-1}\mathbb{E}\Big[\frac{F^s(\widehat{x}^s) - F^s(x_*^s)}{\sigma\alpha m/2} + \big(F^s(\widehat{x}^{s+1}) - F^s(\widehat{x}^s)\big) - 2\big(F(\widehat{x}^{s+1}) - F(\widehat{x}^s)\big)\Big]$$

Above, inequality ① uses Lemma 7.4 and $\Phi^s \geq 0$; inequality ② uses the definition of $\Phi^s$ from Lemma 7.2; inequality ③ uses $\alpha p^2 C(\ell_1 + \ell_2) \leq 1$ and $\alpha p^2 C(\ell_1 \ell_2) \leq 2\sigma$; and equality ④ uses the definition $F^s(y) = F(y) + \sigma \|y - \widehat{x}^s\|^2$.

Finally, rearranging both sides, and using the fact that $\frac{1}{\sigma \alpha m} \leq \frac{1}{6}$, we have

$$\sum_{s=0}^{p-1} \mathbb{E}\left[\left(F^s(\widehat{x}^s) - F^s(x_*^s)\right)\right] \leq 3\mathbb{E}\left[F(\widehat{x}^0) - F(\widehat{x}^p)\right] \ . \qquad \square$$

## 8 Full Method: Final Theorem

We are now ready to state and prove our main convergence theorem for `Natasha1`<sup>full</sup>:

> **Theorem 2.** *Suppose $f(x)$ is $L$-smooth, each $f_i(x)$ is $(\ell_1, \ell_2)$-smooth, $f(x)$ is $\sigma$-nonconvex, and $\sigma \leq \min\{\ell_1, \ell_2, L\}$. If $\frac{\ell_1 \ell_2}{\sigma^2} \leq n$, $p = \Theta\left(\left(\frac{\sigma^2}{\ell_1 \ell_2} n\right)^{1/3}\right)$ and $\Theta\left(\frac{\sigma}{p^2 \ell_1 \ell_2}\right)$, `Natasha1`<sup>full</sup> outputs a point $x^{\mathsf{out}}$ satisfying*
>
> $$\mathbb{E}[\|\mathcal{G}_\eta(x^{\mathsf{out}})\|^2] \leq O\left(\frac{(\ell_1 \ell_2 \sigma)^{1/3} n^{2/3}}{T' n}\right) \cdot (F(x^\varnothing) - F^*)$$
>
> *for $\eta = \frac{1}{\max\{L, 4\sigma\}}$. In other words, to obtain $\mathbb{E}[\|\mathcal{G}_\eta(x^{\mathsf{out}})\|^2] \leq \varepsilon^2$, we need gradient complexity*
>
> $$O\left(n \log \frac{1}{\eta \varepsilon \sigma} + \frac{(\ell_1 \ell_2 \sigma)^{1/3} n^{2/3}}{\varepsilon^2} \cdot (F(x^\varnothing) - F^*)\right) \ .$$

*Remark 8.1.* One can assume $\sigma \leq L$ without loss of generality because any $L$-smooth function is $L$-nonconvex. One can assume $\sigma \leq \ell_2$ without loss of generality because $f(x)$ is $\ell_2$-nonconvex if each $f_i(x)$ is $(\ell_1, \ell_2)$-smooth. Only $\sigma \leq \ell_1$ is a minor requirement for Theorem 2, but if this is not true, one can replace $\ell_1$ with $\sigma$ before applying Theorem 2.

*Remark 8.2.* In Theorem 2 we have assumed $\frac{\ell_1 \ell_2}{\sigma^2} \leq n^2$. If this inequality does not hold, one should apply repeatSVRG instead and it gives faster running time (recall Figure 1a). More specifically, repeatSVRG gives a complexity of

$$\widetilde{O}\left(\frac{n\sigma + n^{1/2}((\ell_1 + \ell_2)\sigma)^{1/2} + n^{3/4}(\ell_1 \ell_2 \sigma^2)^{1/4}}{\varepsilon^2}\right)$$

under a mild assumption of $\sigma \geq \varepsilon^2$ in this more refined $(\ell_1, \ell_2)$-smoothness setting.[13]

*Remark 8.3.* People are also interested in the special case of $\ell_2 = \sigma$ and $\ell_1 = L$. In this case, $L/\sigma$ can be viewed as the "condition number" of the nonconvex problem. In such a case, the gradient complexity of repeatSVRG becomes $\widetilde{O}\left(\frac{\sigma n + \sqrt{nL\sigma}}{\varepsilon^2}\right)$, and that of `Natasha1`<sup>full</sup> becomes $\widetilde{O}\left(n + \frac{(L\sigma^2)^{1/3} n^{2/3}}{\varepsilon^2}\right)$. `Natasha1` is faster when $\frac{L}{\sigma} < n$, and repeatSVRG is faster when $\frac{L}{\sigma} > n$.[14]

*Proof of Theorem 2.* One can verify that our choices of $p$ and $\alpha$ satisfy $p \in [n]$, $\alpha \leq O\left(\frac{\sigma}{p^2 \ell_1 \ell_2}\right)$ and $\alpha \geq \Omega\left(\frac{1}{\sigma m}\right)$, so we can apply Lemma 7.5 and telescope it for the entire algorithm (which has $T'$ full epochs). Use the fact that $\widehat{x}^p$ of the previous epoch equals $\widehat{x}^0$ of the next epoch, we conclude that if we choose a random epoch and a random subepoch $s$, we will have

$$\mathbb{E}[F^s(\widehat{x}^s) - F^s(x^s_*)] \leq \frac{3}{pT'}(F(x^\varnothing) - F^*) \ .$$

---

[13] In V1-V4 of this paper, we forgot to add the additive term $\sqrt{n(\ell_1 + \ell_2)\sigma}$ in the numerator. We have corrected this mistake in version 5.

[14] In this special case, due to our miscalculation in V1-V4, we stated a wrong (and smaller) complexity than $\widetilde{O}\left(\frac{\sigma n + \sqrt{nL\sigma}}{\varepsilon^2}\right)$ for repeatSVRG until September 2018 (version 5). Meanwhile, in May 2018 and September 2018 respectively, Lan and Yang [21] and Chen and Yang [14] discovered different algorithms to achieve this gradient complexity $\widetilde{O}\left(\frac{\sigma n + \sqrt{nL\sigma}}{\varepsilon^2}\right)$. As a result, please be aware that Lan and Yang [21] are the first to state this gradient complexity $\widetilde{O}\left(\frac{\sigma n + \sqrt{nL\sigma}}{\varepsilon^2}\right)$. Stagewise-Katyusha [14] aims to improve the gradient complexity of $\widetilde{O}\left(\frac{\sigma n + n^{3/4}\sqrt{L\sigma}}{\varepsilon^2}\right)$, which was derived from a more general setting in V1-V4 of this work. In the regime of $\frac{L}{\sigma} > n$, however, the current complexity statements of [14, 21] cannot outperform `Natasha1`<sup>full</sup>.

By the $\sigma$-strong convexity of $F^s(\cdot)$, we have $\mathbb{E}[\sigma\|\widehat{x}^s - x^s_*\|^2] \leq \frac{6}{pT'}(F(x^\varnothing) - F^*)$.

Now, $F^s(x) = F(x) + \sigma\|x - \widehat{x}^s\|^2$ satisfies the assumption of $G(x)$ in Lemma 4.1. If we use the SVRG method (see Theorem 2.3) to minimize the convex function $F^s(x)$, we get an output $x^{\mathsf{out}}$ satisfying $F^s(x^{\mathsf{out}}) - F^s(x^s_*) \leq \eta^2 \varepsilon^2 \sigma$ in gradient complexity $O\big((n + \frac{\ell_1 \ell_2}{\sigma^2}) \log \frac{1}{\eta \varepsilon \sigma}\big) \leq O(n \log \frac{L}{\varepsilon \sigma})$.

We can therefore apply Lemma 4.1 and conclude that this output $x^{\mathsf{out}}$ satisfies

$$\mathbb{E}[\|\mathcal{G}_\eta(x^{\mathsf{out}})\|^2] \leq O\Big(\frac{\sigma}{pT'}\Big) \cdot (F(x^\varnothing) - F^*) = O\Big(\frac{(\ell_1 \ell_2 \sigma)^{1/3} n^{2/3}}{T'n}\Big) \cdot (F(x^\varnothing) - F^*) \ .$$

In other words, we obtain $\mathbb{E}[\|\mathcal{G}_\eta(x^{\mathsf{out}})\|^2] \leq \varepsilon^2$ using

$$T'n = O\Big(n + \frac{(\ell_1 \ell_2 \sigma)^{1/3} n^{2/3}}{\varepsilon^2} \cdot (F(x^\varnothing) - F^*)\Big)$$

computations of the stochastic gradients. Here, the additive term $n$ is because $T' \geq 1$.

Finally, adding this with $O(n \log \frac{L}{\varepsilon \sigma})$, the gradient complexity for the application of SVRG in the last line of `Natasha1`$^{\mathsf{full}}$, we finish the proof of the total gradient complexity. □

## Acknowledgements

We thank Yuanzhi Li for enlightening conversations, Weiran Wang for correcting typos, and Zaiyi Chen, Tianbao Yang and Guanghui Lan for correcting a mistake when stating the complexity of repeatSVRG. This paper was partially supported by a Microsoft Research Award, no. 0518584, and an NSF grant, no. CCF-1412958. Any opinions, findings and conclusions or recommendations expressed in this material are those of the authors and do not reflect the views of the National Science Foundation and Microsoft Corporation.

# APPENDIX

## A  Other Related Works

Methods based on variance-reduced stochastic gradients were first introduced for convex optimization. The first such method is SAG by Schmidt et al [29]. The two most popular choices for gradient estimators are the SVRG-like one we adopted in this paper (independently introduced by [20, 32], and the SAGA-like one introduced by [15]. In nearly all applications, the results proven for SVRG-like estimators and SAGA-like estimators are simply exchangeable (therefore, the results of this paper naturally generalize to SAGA-like estimators).

The first "non-convex use" of variance reduction is by Shalev-Shwartz [30] who assumes that each $f_i(x)$ is non-convex but their average $f(x)$ is still convex. This result has been slightly improved to several more refined settings [12]. The first truly non-convex use of variance reduction (i.e., for $f(x)$ being also non-convex) is independently by [6] and [28]. First-order methods only find stationary points (unless there is extra assumption on the randomness of the data), and converge no faster than $1/\varepsilon^2$.

When the second-order Hessian information is used, one can (1) find local minima instead of stationary points, and (2) improve the $1/\varepsilon^2$ rate to $1/\varepsilon^{1.5}$. The first such result is by cubic-regularized Newton's method [26]; however, its per-iteration complexity is very high. Very recently, two independent groups of authors tackled this problem from a somewhat similar viewpoint [1, 13]: if the computation of Hessian-vector multiplications (i.e., $(\nabla^2 f_i(x))v$) is on the same order of the computation of gradients $\nabla f_i(x)$,[15] then one can obtain a $(\varepsilon, \sqrt{\varepsilon})$-approximate local minimum in gradient complexity $\widetilde{O}\big(\frac{n}{\varepsilon^{1.5}} + \frac{n^{3/4}}{\varepsilon^{1.75}}\big)$, if we use big-$O$ to also hide dependencies on the smoothness parameters. Although Carmon et al. [13] only stated a complexity of $\widetilde{O}\big(\frac{n}{\varepsilon^{1.75}}\big)$ in the non-stochastic setting, their result generalizes to our stated complexity in the stochastic setting. As we argue in Appendix B, both these methods reduce the problem of finding $(\varepsilon, \sqrt{\varepsilon})$-approximate local minima to that of finding $\varepsilon$-approximate stationary points in functions of $\sqrt{\varepsilon}$-bounded nonconvexity.

Other related papers include Lee et al. [22] who showed that gradient descent, starting from a random point, almost surely converges to a local minimum if the function is "strict-saddle". The rate of convergence required is somewhat unknown.

**Online Stochastic Setting.** Ge et al. [18] where the authors showed that a noise-injected version of SGD converges to local minima instead of critical points, as long as the underlying function is "strict-saddle." Their theoretical running time is a large polynomial in the dimension.

# B  From Stationary Points to Local Minima

Recently, researchers have shown that the general problem of finding $(\varepsilon, \delta)$-approximate local minima, under mild conditions, reduces to (repeatedly) finding $\varepsilon$-approximate stationary points for an $O(\delta)$-nonconvex function [1, 3, 5, 13]. The authors of [3, 5] call this approach "swing by saddle points," which is different from the classical approach to escape from saddle points. We sketch the details here for the sake of completeness, in the special case of $\psi(x) \equiv 0$.[16]

We say that a point $x$ is $(\varepsilon, \delta)$-approximate local minimum, if $\|\nabla f(x)\| \leq \varepsilon$ and $\nabla^2 f(x) \succeq -\delta \mathbf{I}$.

Carmon et al. [13] showed that an $(\varepsilon, \delta)$-approximate minimum for the general problem (1.1) can be solved via the following iterative procedure. In every iteration at point $x_t$, detect whether the smallest eigenvalue of $\nabla^2 f(x_t)$ is below $-\delta$:

- if yes, find the smallest eigenvector of $\nabla^2 f(x_t)$ approximately and move in this direction. (One can use for instance the shift-and-invert method [17], momentum method [10], Oja's method [9, 10], or the slower power method [31].)

- if no, define $f_t(x) = f(x) + L\big(\max\big\{0, \|x - x_t\| - \frac{\delta}{L_2}\big\}\big)^2$ where $L_2$ is the second-order smoothness of $f(x)$ and $f_t(x)$ can be proven as $5L$-smooth and of $3\delta$-bounded nonconvexity; we then find an $\varepsilon$-approximate stationary point of $f_t(x)$ and move there.

**The Trade-Off on $\delta$.** The final running time of the above algorithm depends on the maximum between (1) the eigenvector computation and (2) the stationary-point computation. The larger $\delta$ is, the faster (1) becomes and the slower (2) becomes; the smaller $\delta$ is, the faster (2) becomes and the slower (1) becomes.

In the offline setting, as argued in [1, 13], the optimum trade-off is $\delta = \sqrt{L_2 \varepsilon}$. This confirms that in optimization with $\sigma$-bounded nonconvexity, the parameter $\sigma$ can be much smaller than $L$.

---
[15]A lot of interesting problems satisfy this property, including training neural nets.
[16]This reduction was first explicitly given in [13] but only implicitly in [1].